\documentclass[preprint]{iucr}              

                   \def\href#1{\relax}\let\foo\caption
                   
\usepackage{url}

\usepackage{breakurl}
                   
\ifPDF
  \RequirePackage{hyperref}
  \PassOptionsToPackage{pdftex,bookmarksopen,bookmarksnumbered,breaklinks}{hyperref}
\fi
\let\caption\foo

     \paperprodcode{a000000}      
     \paperref{xx9999}            
     \papertype{IU}               

     \paperlang{english}          
     \journalyr{2003}
     \journalreceived{\relax}
     \journalaccepted{\relax}
     \journalonline{\relax}
\usepackage{amssymb}
\usepackage{amsthm}
\usepackage{chemformula}
\usepackage[version=3]{mhchem} 
\usepackage{physics}
\usepackage{verbatim}

\theoremstyle{definition}
\newtheorem{thm}{Theorem}
\newtheorem{dfn}[thm]{Definition}
\newtheorem{lem}[thm]{Lemma}

\newcommand{\R}{\mathbb{R}}
\newcommand{\Z}{\mathbb{Z}}

\newcommand{\SL}{\mathrm{SL}}
\newcommand{\PDD}{\mathrm{PDD}}
\newcommand{\PDA}{\mathrm{PDA}}

\newcommand{\ADA}{\mathrm{ADA}}

\newcommand{\CIA}{\mathrm{CIA}}
\newcommand{\aCIA}{\overline{\CIA}}
\newcommand{\EMD}{\mathrm{EMD}}
\newcommand{\RMS}{\mathrm{RMS}}

\newcommand{\PPC}{\mathrm{PPC}}

\newcommand{\ep}{\varepsilon}

\newcommand{\al}{\alpha}

\newcommand{\La}{\Lambda}
\newcommand{\ga}{\gamma}
\newcommand{\be}{\beta}

\newcommand{\vol}{\mathrm{vol}}

\newcommand{\angstrom}{\textup{\AA}}

\newcolumntype{R}[1]{>{\raggedleft\let\newline\\\arraybackslash\hspace{0pt}}m{#1}}
\newcolumntype{L}[1]{>{\raggedright\let\newline\\\arraybackslash\hspace{0pt}}m{#1}}
\newcolumntype{C}[1]{>{\centering\let\newline\\\arraybackslash\hspace{0pt}}m{#1}}

\graphicspath{{images/}}

\begin{document}                  



\title{Continuous invariant-based asymmetries of periodic crystals quantify deviations from higher symmetry}
\shorttitle{Continuous invariant-based asymmetries of crystals}


     \author[a]{Surya}{Majumder}
     \author[b]{Daniel}{Widdowson}
     \author[a]{Yury}{Elkin}
     \author[a,c]{Olga}{Anosova}
     \author[b]{Andrew I.}{Cooper}
     \author[d]{Graeme M.}{Day}
     \cauthor[a,b,c]{Vitaliy A.}{Kurlin}{vkurlin@liv.ac.uk}{}\aufn{}
\aff[a]{Computer Science department, University of Liverpool, Liverpool, L69 3BX, UK}\aff[b]{Materials Innovation Factory, University of Liverpool, Liverpool, L7 3NY, UK}     \aff[c]{National Institute for Theory and Mathematics in Biology, Chicago, US}
\aff[d]{School of Chemistry and Chemical Engineering, University of Southampton,
Southampton, SO17 1NX, UK}
    

\shortauthor{Majumder et al}







\maketitle                        

\begin{synopsis}
The new continuous invariant-based asymmetry quantifies a deviation of any periodic crystal from its closest higher symmetry neighbour where all molecules are geometrically equivalent.
\end{synopsis}

\vspace*{-10mm}

\begin{abstract}
Ideal symmetry is known to break down under almost any noise.
One measure of asymmetry in a periodic crystal is the relative multiplicity $Z'$ of
geometrically non-equivalent units.
However, $Z'$ discontinuously changes under almost any displacement of atoms, which
can arbitrarily scale up a primitive cell.
This discontinuity was recently resolved by a hierarchy of invariant descriptors that
continuously change under all small perturbations.
We introduce a Continuous Invariant-based Asymmetry (CIA) to quantify (in physically
meaningful Angstroms) the deviation of a periodic crystal from a higher symmetry
form.
Our experiments on several Crystal Structure Prediction datasets show that about a
half of simulated crystals have high values of $\CIA$, while all experimental
structures in these datasets have $\CIA=0$.
On another hand, many crystals with high values $Z'$ in the Cambridge Structural
Database (CSD) turned out to be close to more symmetric forms with $Z'\leq 1$ due to
low values of $\CIA$s.
\end{abstract}


\section{Introduction: motivations for a new continuous asymmetry of crystals}
\label{sec:intro}

Many periodic crystals are highly symmetric, because a globally stable structure is
usually formed by a few energetically favourable interactions, bonds, molecules, or
formula units, which are repeated in $\R^3$ by symmetries \cite{lax2001symmetry}.
Though we were motivated by molecular crystals, our invariant-based approach to
asymmetry extends to all non-molecular crystals and periodic sets in any Euclidean
space $\R^n$.
\smallskip

While molecular crystals can contain many molecules in primitive unit cells, they are
often obtained from a smaller number of molecules by \emph{symmetry operations}
preserving the whole crystal in $\R^3$ \cite{Zprime}.
For a non-molecular crystal, the chemical analogue of a molecule is a \emph{formula
unit} that is an electronically neutral group of atoms or ions, embedded in $\R^3$
and representing their relative numbers in a given compound, reduced to the smallest
integer numbers.
For example, table salt has the empirical formula \ce{NaCl} with a formula unit consisting of two ions \ce{Na}$^+$ and \ce{Cl}$^-$.
However, this pair of ions can be chosen in many geometrically different ways, because ionic bonds in \ce{NaCl} do not define a bounded object, such as a molecule.
Hence, non-molecular crystals should be unambiguously split into disjoint \emph{geometric blocks}, for example, single ions, or metal blocks and organic linkers in a metal organic framework.
\smallskip

In this paper, a \emph{crystal} $S$ means a periodic crystal, while $Z$ can be
non-integer for disordered or aperiodic crystals \cite{senechal1996quasicrystals}.
The \emph{multiplicity} $Z$ usually denotes the number of formula units in a
primitive unit cell.
The \emph{relative multiplicity} \emph{Z prime}
was often defined as $Z(S)$ divided by the number of independent general positions
\cite{steed2015packing}, which makes sense, if $S$ consists of chemically equivalent molecules.
For crystals with chemically different molecules (called co-crystals), \cite{van2000structure} used another notation $Z''$ for the number of crystallographically non-equivalent molecules.
To cover non-molecular crystals, we define $Z'$ below for any periodic point set $S\subset\R^n$ with a given splitting into geometric blocks.

\begin{dfn}[relative multiplicity $Z'$]
\label{dfn:Zprime}
An \emph{asymmetric unit} is a minimal, closed, and simply connected subset $A$ of a
unit cell of $S\subset\R^b$, whose images under all symmetry operations of $S$ tile
the whole space $\R^n$.
Let $S\cap A$ split into geometric blocks $B_1,\dots,B_G$, which should be chemically
different molecules, ions, or other disjoint blocks for crystals in $\R^3$.
Let $B_i$ have $n_i$ symmetry operations (including the identity) that preserve both
$S$ and $B_i$, $i=1,\dots,G$.
The \emph{relative multiplicity} is $Z'(S)=\sum\limits_{i=1}^{G} \dfrac{1}{n_i}$.
\end{dfn}

Geometric blocks $B_i,B_j$ of $S\subset\R^n$ are called \emph{rigidly equivalent} if
there is a rigid motion of $\R^n$ that maps $S$ to $S$ and $B_i$ to $B_j$.
If all molecules of a crystal $S$ are rigidly equivalent, an asymmetric unit $A$ of
$S$ contains one molecule $B$, so $Z'(S)=\dfrac{1}{n}=\dfrac{Z}{N}$, where $n$ is the
number of symmetry operations preserving both $S$ and $B$, while $N=nZ$ is the number
of symmetry operations preserving $S$ and the motif $S\cap U$, which can permute
molecules within a primitive cell $U$ of $S$.
If $S\cap A$ consists of two non-rigidly equivalent molecules in 2-fold positions,
then
$Z'(S)=\dfrac{1}{2}+\dfrac{1}{2}=1$.
The crystal \ce{NaCl} has ions \ce{Na}$^+$ and \ce{Cl}$^-$ with point groups of order
48, so $Z'(\ce{NaCl})=\dfrac{1}{48}+\dfrac{1}{48}=\dfrac{1}{24}$.
\smallskip
 
In about 90\% of entries in the Cambridge Structural Database (CSD), an asymmetric
unit includes only one molecule, so $Z'\leq 1$ \cite{anderson2006when}.
However, the CSD has many crystals with high $Z'$ \cite{brock2016high}, e.g. OGUROZ
has $Z'=56$.
\smallskip

\begin{figure}
\centering
\includegraphics[width=\textwidth]{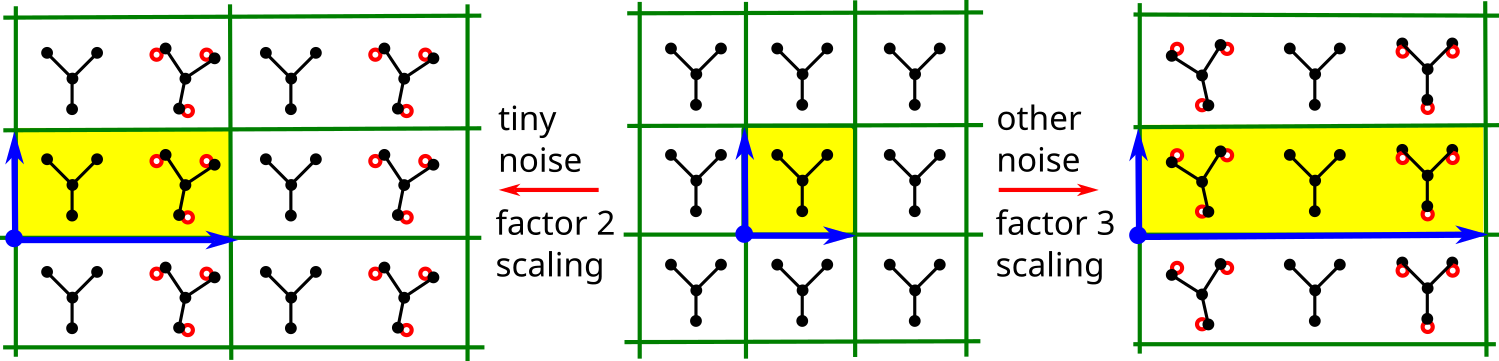}
\caption{Almost any noise arbitrarily scales up a primitive yellow cell and
discontinuously changes the relative multiplicity $Z'$ of molecules, which are
represented by black $Y$ graphs whose terminal vertices have initial positions shown
by red circles.}
\label{fig:noise_scales_cells_Z'}
\end{figure}

Crystal Structure Prediction (CSP) often starts with simulating $Z'=1$ crystals for
the most frequent space groups, but a final energy relaxation can produce structures
with $Z'$ values up to 36 \cite{pulido2017functional}.
More importantly, almost any displacement of atoms or a whole rigid molecule
discontinuously changes the size of a primitive (or reduced) cell and hence
arbitrarily scales up $Z'$.
Fig.~\ref{fig:noise_scales_cells_Z'} shows nearly identical structures with
$Z'=1,2,3$ and similarly applies to any periodic crystal.
\smallskip

Ignoring any noise up to a small threshold $\ep$ only shifts the problem from $0$ to
another number without guarantees of a continuous change.
This \emph{sorites} paradox (when a heap of sand stops being a heap while grains are
removed one by one) has been known since ancient times \cite{sorites}.
Its rigorous solution requires an \emph{invariant} that is preserved by any rigid
transformation and continuously changes under perturbations of atoms.
\smallskip

While a full hierarchy of such invariants for periodic crystals from the fastest to
complete is being finalised by \cite{anosova2025geometric,widdowson2025higher}, our
continuous asymmetry will be based on the fast invariant $\PDD$ (Pointwise Distance
Distribution), which distinguishes all non-duplicate crystals in the world's largest
databases within two hours on a modest desktop \cite{widdowson2022resolving}.

\section{Generically complete and continuous isometry invariants of crystals}
\label{sec:invariants}

This section recalls isometry invariants, which will be used to define a continuous
invariant-based asymmetry in section~\ref{sec:asymmetry}.
Any linear basis $\vb*{v_1},\dots,\vb*{v_n}$ of $\R^n$ generates \\
the \emph{lattice} $\La=\{c_1\vb*{v_1}+c_2\vb*{v_2}+\dots+c_n\vb*{v_n} \mid
c_1,\dots,c_n\in\Z\}\subset\R^n$ and \\
the \emph{unit cell} $U=\{x_1\vb*{v_1}+x_2\vb*{v_2}+\dots+x_n\vb*{v_n} \mid 0\leq
x_1,\dots,x_n<1\}\subset\R^n$.
\smallskip

For any finite motif $M\subset U$ of atoms (considered zero-sized points) in the unit
cell, a periodic crystal is defined as the infinite set $S=\La+M=\{\vb*{v}+p \mid
\vb*{v}\in\La, p\in M\}$ or a finite union $\cup_{p\in M} (\La+p)$ of shifted
lattices with origins at all points of $M$.
\smallskip
 
The definition above is widely used for representing crystals in Crystallographic
Information Files (CIFs), but is highly ambiguous in the sense that infinitely many
pairs (basis, motif) represent the same crystal $S$.
This ambiguity motivated us to distinguish between a crystal $S$ and its
\emph{structure}, defined as the equivalence class of all periodic sets of atoms that
are represented by different CIFs but can be exactly matched with each other by rigid
motion, see Definition~6 in \cite{anosova2024importance}.
\smallskip

Any canonical (standard or conventional) choice of a cell for a periodic crystal is
discontinuous under almost any noise, as in Fig.~\ref{fig:noise_scales_cells_Z'},
which was experimentally demonstrated already in 1965, see p.80 in
\cite{lawton1965reduced}.
The new definition of a \emph{crystal structure} as a rigid class (consisting of all
crystals that can be exactly matched under rigid motion) has become practical due to
the hierarchy of invariants that uniquely identify any crystal structure independent
of its initial representation.
\smallskip

Definition~\ref{dfn:PDD} introduces the invariant $\PDD$ for any periodic set of
points in $\R^n$, which can be all atomic centres of a crystal in $\R^3$, or other
points defined by a crystal, for example, atoms of one specific type, or molecular
centres, which form a periodic set.

\begin{dfn}[Pointwise Distance Distribution $\PDD$]
\label{dfn:PDD}
Let $S \subset \mathbb{R}^n$ be a periodic point set with a motif $M = \{p_1,
p_2,\dots, p_m\}$.
Fix an integer $k\ge 1$. 
For every point $p_i \in M$, let $d_1(p)\le\cdots\le d_k(p)$ be the distances from
$p$ to its $k$ nearest neighbours within the full infinite set $S$ not restricted to
any cell.
The matrix $D(S;k)$ has $m$ rows consisting of the distances $d_1(p_i),\dots,
d_k(p_i)$ for $i=1,\dots,m$.
If any $l\ge 1$ rows are identical to each other, we collapse them into a single row
and assign the weight $\dfrac{l}{m}$ to this row.
The resulting matrix of $k$ columns and a maximum of $m$ rows with weights, in the
extra (say, 0-th) column, is called the Pointwise Distance Distribution $\PDD(S;k)$.
\end{dfn}

The columns of the matrix $\PDD(S;k)$ are ordered because each row consists of
increasing values of distances to neighbours but without their indices.
So $\PDD(S;k)$ importantly differs from the matrix of pairwise distances between $m$
points in the motif $M$, also because neighbours are not restricted to any (extended)
cell of $S$.
\smallskip

Since many crystals consist of indistinguishable atoms, we consider all points of $S$ unordered.
Then $\PDD(S;k)$ has unordered rows and can be interpreted as a discrete distribution
of rows (or unordered points in $\R^k$) with probabilities equal to the weights
assigned to rows.
The Pair Distribution Function is obtained from a single collection of all interatomic distances (usually normalised by frequencies and then smoothed) and hence
is naturally weaker than $\PDD(S;k)$, which splits distances per point and avoids
losing information under smoothing, see the discussion at the end of section 3 in
\cite{widdowson2022resolving}.
This probabilistic interpretation allows one to compare $\PDD$s by many distance
metrics on discrete distributions.
We usually use the simplest metric called Earth Mover's Distance (EMD), which was
adapted for comparing chemical compositions \cite{hargreaves2020earth}.
Theorem 4.2 in \cite{widdowson2026pointwise} proved that $\PDD(S;k)$ continuously
changes in $\EMD$ under perturbations, including those that arbitrarily scale up a
minimal cell as in Fig.~\ref{fig:noise_scales_cells_Z'}.
\smallskip

The most important result about the $\PDD$ is its generic completeness: Theorem~5.8
in \cite{widdowson2026pointwise} proved that $\PDD(S;k)$ with a lattice of $S$ and
the number $m$ of points in a motif of $S$ suffice to reconstruct any generic
periodic point set $S\subset\R^n$, uniquely under isometry, for a large enough $k$
with an explicit upper bound.
In other words, $\PDD(S;k)$ with a few extra invariants provably distinguishes all
crystals, possibly except singular examples that form a subspace of measure 0 within
the continuous space of all periodic crystals.
In practice, $\PDD(S;k)$ distinguished all non-duplicate crystals in the world's
major databases within two hours on a modest desktop, see Table~3 in
\cite{widdowson2026pointwise}.
Theorem 3.7 in \cite{widdowson2026pointwise} proved that, as $k\to+\infty$, the
distances in each row of $\PDD(S;k)$ asymptotically approach $\PPC(S)\sqrt[n]{k}$,
where the Point Packing Coefficient $\PPC(S)$ is inversely proportional to the point
density, as defined below.
This fact motivated us to subtract this asymptotic curve from $\PDD(S;k)$ to
neutralise the influence of density.

\begin{dfn}[invariants $\PPC(S)$ and $\PDA(S;k)$]
\label{dfn:PDA}
Let $S\subset\R^n$ be a periodic set with $m$ points in a unit cell $U$ of $S$. 
The \emph{Point Packing Coefficient} is $\PPC(S) = \sqrt[n]{\dfrac{\vol(U)}{mV_n}}$,
where $\vol(U)$ is the volume of $U$, and $V_n$ is the volume of the unit ball in
$\mathbb{R}^n$, e.g. $V_3=\frac{4}{3}\pi$.
The \emph{Pointwise Deviation from Asymptotic} is the matrix $\PDA(S;k)$ obtained
from $\PDD(S;k)$ by subtracting $\PPC(S)\sqrt[n]{j}$ from every distance in columns
$j=1,\dots,k$.
\end{dfn}

Another advantage of $\PDA(S;k)$ vs original $\PDD(S;k)$ is the experimental
convergence to 0 of the $k$-th values from the last column of $\PDA(S;k)$ as
$k\to+\infty$, see Fig.~4 in \cite{widdowson2025geographic}.
This convergence to 0 was formally proved for any cubic lattice $\Z^n$ in Example~SM3.1 from \cite{widdowson2026pointwise}.
Then there is no need to substantially increase the number $k$ of neighbours, because more distant neighbours bring smaller contributions.
We consider $k$ not as a parameter that seriously affects $\PDA(S;k)$, but as a degree of approximation like the number of decimal places on a calculator.
The vector $\ADA(S;k)$ of column averages in $\PDA(S;k)$ for $k=100$ atomic neighbours was sufficient to distinguish all non-duplicate crystals in the CSD
\cite{widdowson2024continuous}.
The components of this \emph{Average Deviation from Asymptotic} vector $\ADA(S;k)$ can be used as analytic coordinates on geographic-style maps of any materials
database.
Such maps were first developed for 2D lattices by
\cite{bright2023geographic,bright2023continuous,kurlin2024mathematics}.

\section{A continuous invariant-based asymmetry (CIA) of periodic crystals}
\label{sec:asymmetry}

The discontinuity of $Z'$ from Definition~\ref{dfn:Zprime} under almost any
perturbation has been known for 30+ years.
The quote ``two fairly unsymmetrical objects can be combined into a less
unsymmetrical structural dimer'' from \cite{wilson1993space} means that a crystal
with $Z'=2$ can be geometrically close to a more symmetric crystal with $Z'=1$.
\smallskip

This section first defines the Earth Mover's Distance ($\EMD$) between geometric
blocks within a periodic point set $S\subset\R^n$ by using the isometry invariant
$\PDA(S;k)$ from Definition~\ref{dfn:PDA}.
The continuous invariant-based asymmetry of $S$ will be defined through $\EMD$s
between all blocks in an asymmetric unit of $S$.
The $\EMD$ needs a ground distance between vectors $\vb*{b}=(b_1,\dots,b_k)$ and
$\vb*{c}=(c_1,\dots,c_k)$ in $\R^k$, such as rows of $\PDA(S;k)$.
The simplest choices are the \emph{Chebyshev} distance
$d_\infty(\vb*{b},\vb*{c})=\max\limits_{1\leq i\leq k}|b_i-c_i|$ and the \emph{Root
Mean Square} ($\RMS$)
$d(\vb*{b},\vb*{c})=\sqrt{\dfrac{1}{k}\sum\limits_{i=1}^k (b_i-c_i)^2}$.
\smallskip

These distances respect the continuity under perturbations as follows.
If any $b_i,c_i$ are perturbed up to $\ep$, then
$|b_i-c_i|\leq 2\ep$ for $i=1,\dots,k$, and
both $d_\infty(\vb*{b},\vb*{c})\leq 2\ep$, 
$d(\vb*{b},\vb*{c})\leq 2\ep$. 
We usually write $d$ without a subscript for brevity.
If $d_\infty$ is used in computations, all relevant distances and asymmetry will have
the subscript $\infty$.
\smallskip

For any periodic set in $\R^n$, Definition~\ref{dfn:EMD} introduces a distance
between geometric blocks $B,C$ (considered as finite sets of points), which are
molecules, ions, or other well-defined disjoint subsets for crystals in $\R^3$.
This distance measures how the positions of $B,C$ differ within a common periodic set
$S$ containing both $B,C$.
If $B,C$ can be exactly matched by a rigid motion of $\R^n$ preserving $S$, then this
distance is 0.
In all real examples, any deviation from symmetry should be positive because of
noise.
\smallskip

Though the $\EMD$ makes sense for distributions of different sizes, our experiments
on crystals will use the $\EMD$ only for geometric blocks that are chemically
identical molecules.
More generally, we assume that every point in a periodic set $S\subset\R^n$ has a
categorical label, which is an analogue of an atomic type, such as \ce{Na}$^+$ and \ce{Cl}$^-$.
\smallskip

Briefly, the $\EMD$ optimally splits and transports objects from one distribution to
another by minimising the overall cost based on a ground distance between objects.
If we need to guarantee matching of points only with the same label (atomic type for
crystals), the ground distance can be adjusted by taking the maximum of $d_\infty$ or
$d=\RMS$ with a discrete metric that is infinite between points of different labels.

\begin{dfn}[Earth Mover's Distance $\EMD$ between geometric blocks]
\label{dfn:EMD}
Let $S\subset\R^n$ be a periodic set of labeled points with an asymmetric unit $A$.
Let $B,C\subset S\cap A$ be geometric blocks (finite sets) that have $m(B),m(C)$
points of weights $\dfrac{1}{m(B)},\dfrac{1}{m(C)}$, respectively.
For $i=1,\dots,m(B)$ and $j=1,\dots,m(C)$, let $R_i(B),R_j(C)$ be the rows of $i$-th
and $j$-th points in $B,C$, respectively.
The distance below is independent of point ordering.
The \emph{Earth Mover's Distance}  
$\EMD(B,C)=\sum\limits_{i=1}^{m(B)} \sum\limits_{j=1}^{m(C)} f_{ij} d(R_i(B),R_j(C))$
is minimised over variable parameters $f_{ij}\in[0,1]$ subject to 
$\sum\limits_{j=1}^{m} f_{ij}=\dfrac{1}{m(B)}$ and 
$\sum\limits_{j=1}^{m} f_{ij}=\dfrac{1}{m(C)}$ for all $i=1,\dots,m(B)$ and
$j=1,\dots,m(C)$, respectively.
\end{dfn}

The distance $\EMD(B,C)$ measures the minimum perturbation of the rows of the geometric blocks $B,C$ in $\PDA(S;k)$ to match (distance-based invariants of) $B$ and $C$ within the ambient periodic set $S$.
This perturbation matching $B$ and $C$ reduces the number of geometrically non-equivalent blocks and hence makes $S$ more symmetric.
\smallskip

If an asymmetric unit $A$ of $S$ has only one geometric block $B$, then $S$ has no asymmetry because all blocks in $S$ are images of $B$ under symmetry operations of $S$.
If $A$ has only two blocks $B,C$, then $\EMD(B,C)$ is considered the asymmetry of $S$.
In more general cases, Definition~\ref{dfn:CIA} introduces the continuous asymmetry below.

\begin{dfn}[Continuous Invariant-based Asymmetry $\CIA(S)$]
\label{dfn:CIA}
Let a periodic set $S\subset\R^n$ with labeled points have geometric blocks
$B_1,\dots,B_G$ its asymmetric unit.
Set $d_{i}=\max\limits_{j=1,\dots,G} \EMD(B_i,B_j)$ for $i=1,\dots,G$.
The \emph{Continuous Invariant-based Asymmetry} is $\CIA(S) =
\min\limits_{i=1,\dots,G} d_{i}$.
The `average' version is $\aCIA(S)=
\dfrac{1}{G}\sum\limits_{i=1}^{G}  d_{i}$.
\end{dfn}

The matrix of distances $\EMD(B_i,B_j)$ describes the relative positions of $G$
blocks within an asymmetric unit of $S$ in terms of their distances to atomic
neighbours within the full $S$.
For $i=1,\dots,G$, the distance $d_i$ measures how far $B_i$ is from all other blocksThe standard (min-max) formula of $\CIA(S)$ means that the optimal $i$-th block $B_i$
serves as a centre minimising its distance $\EMD(B_i,B_j)$ to the farthest block
$B_j$, while $\aCIA(S)$ averages maximum deviations $d_i$ from all blocks considered
as centres.
The default notation $\CIA(S)$ uses $\EMD$ based on the ground distance $d=\RMS$
between rows of $\PDA(S;k)$ with $k=100$.
For the Chebyshev distance $d_\infty$, we keep the subscript $\infty$ in the
notations $\EMD_\infty$, $\CIA_\infty$, and $\aCIA_\infty$.

\begin{lem}[invariance of $\CIA$s]
\label{lem:CIA_invariance}
All $\CIA$s in Definition~\ref{dfn:CIA} are invariant (remain unchanged) under any
isometry and changes of a unit cell of a periodic point $S\subset\R^n$.
\end{lem} 

Lemma~\ref{lem:CIA_invariance} and all further results below are proved in
appendix~\ref{sec:proofs}.

\begin{lem}[inequalities for $\CIA$s]
\label{lem:CIA_inequalities}
In the notations of Definition~\ref{dfn:CIA}, 
$\CIA\leq\CIA_\infty$, $\aCIA\leq\aCIA_\infty$, and
$\CIA\leq\aCIA\leq 2\CIA$ hold for any periodic point set $S\subset\R^n$.
\end{lem} 

Since Definition~\ref{dfn:CIA} is based on the invariant $\PDA(S;k)$, the full
notation should be $\CIA(\PDA(S;k))$, where $\PDA(S;k)$ can be replaced with another
``pointwise'' invariant, such as the higher-order $\PDA^{(h)}$
\cite{widdowson2025higher} or complete isoset \cite{anosova2025recognition}.
In this paper, we use only $\PDA(S;100)$ and write $\CIA(S)$ for brevity.
Theorem~\ref{thm:CIA} justifies the continuity of the asymmetry $\CIA(S)$ under all
small perturbations of points, including those that arbitrarily scale up an initial
cell of $S$.

\begin{thm}[continuity of $\CIA$ under perturbations]
\label{thm:CIA}
Let $S\subset\R^n$ be a periodic point set and $r(S)$ denote the minimum
half-distance between any points of $S$.
For any $0<\ep<r(S)$, let a periodic set $Q\subset\R^n$ be obtained by perturbing
every point of $S$ up to Euclidean distance $\ep$.
Then the $\CIA$s based on the invariant $\PDA(S;k)$ for any $k$ in
Definition~\ref{dfn:CIA} satisfy
$|\CIA(S)-\CIA(Q)|\leq 4\ep$ and
$|\aCIA(S)-\aCIA(Q)|\leq 4\ep$.  
\end{thm}


\section{Fast detection of asymmetric crystals in large simulated CSP datasets}
\label{sec:simulated}

This section visualises several versions of $\CIA$ for 50+ thousand simulated crystals from four CSP datasets reported in \cite{pulido2017functional}.
At that time, the synthesised crystals predicted by these CSPs substantially extended the small population of nanoporous crystals in the CSD.
However, these predictions took more than 12 weeks on a supercomputer, also due to predictions of properties, such as gas capture.
\smallskip

\begin{figure}
    \centering
    \includegraphics[height=32mm]{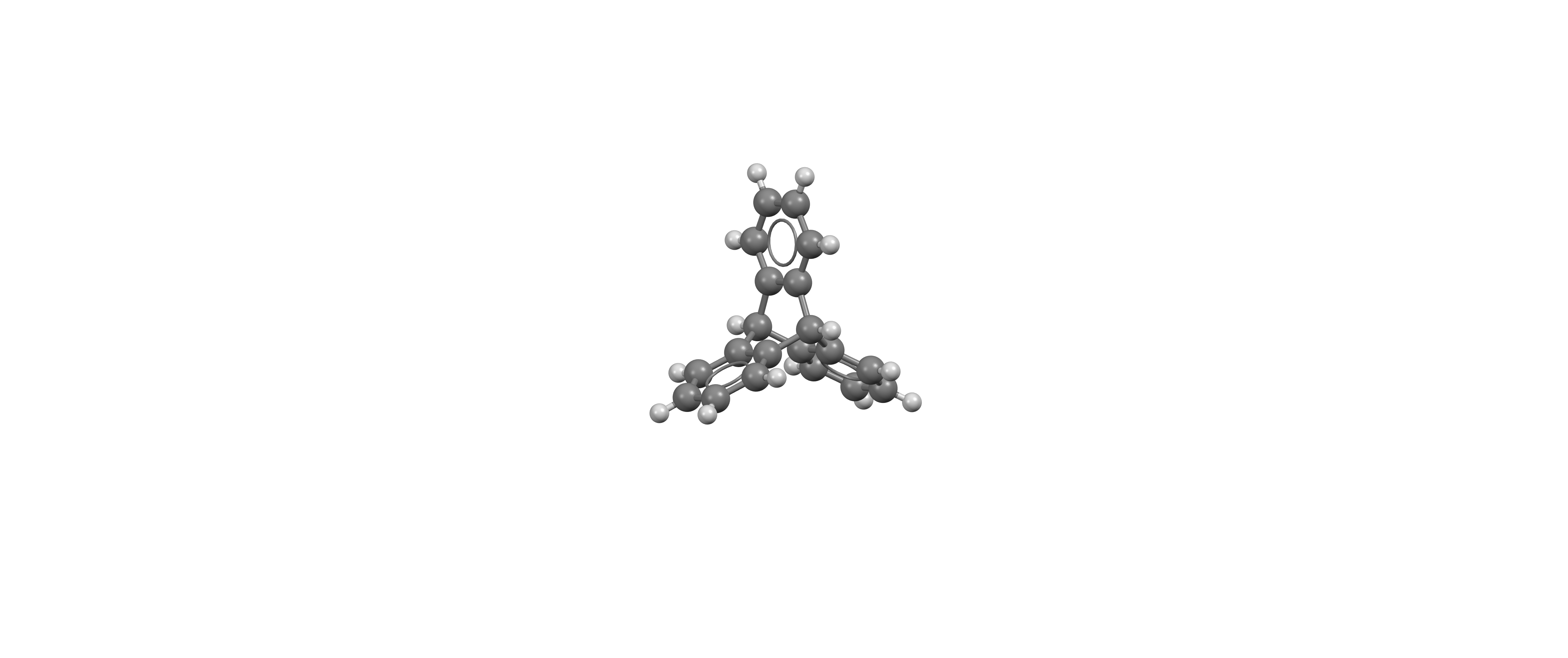}
    \includegraphics[height=32mm]{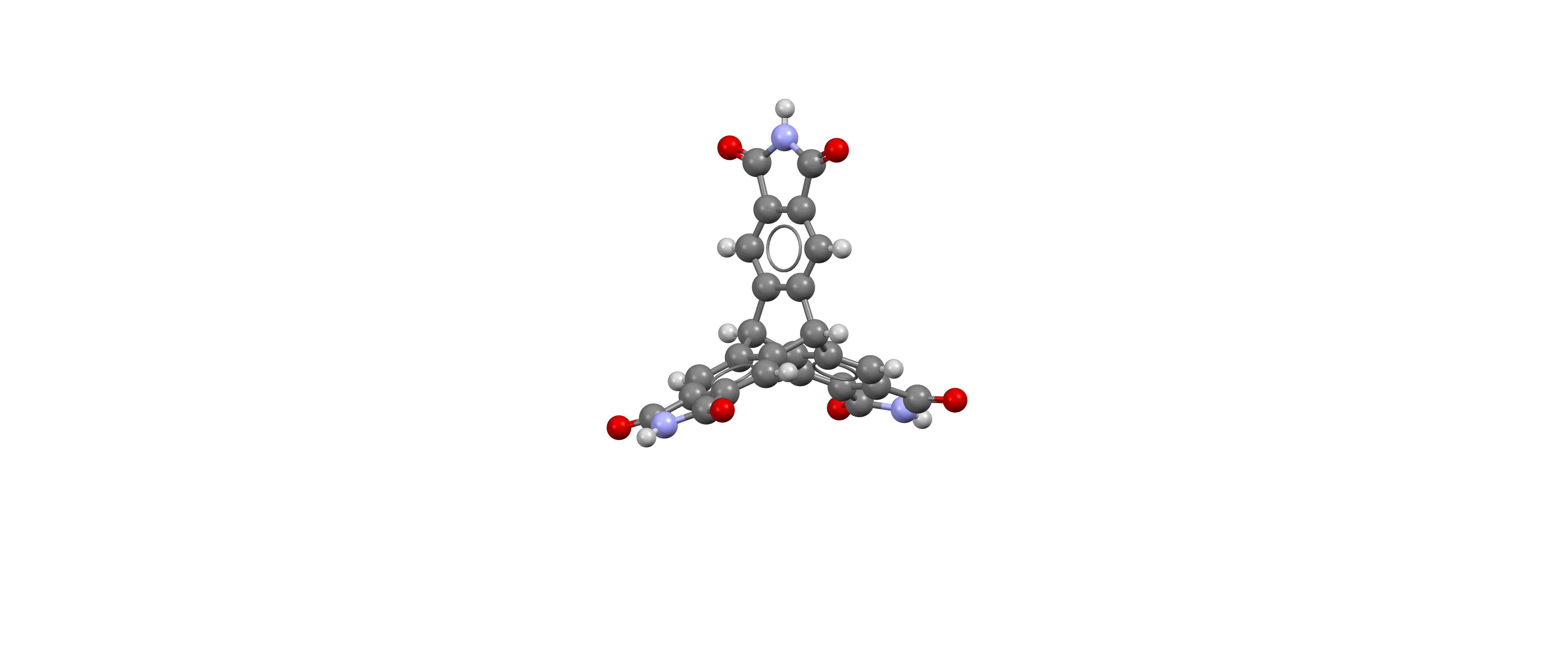}
    \includegraphics[height=32mm]{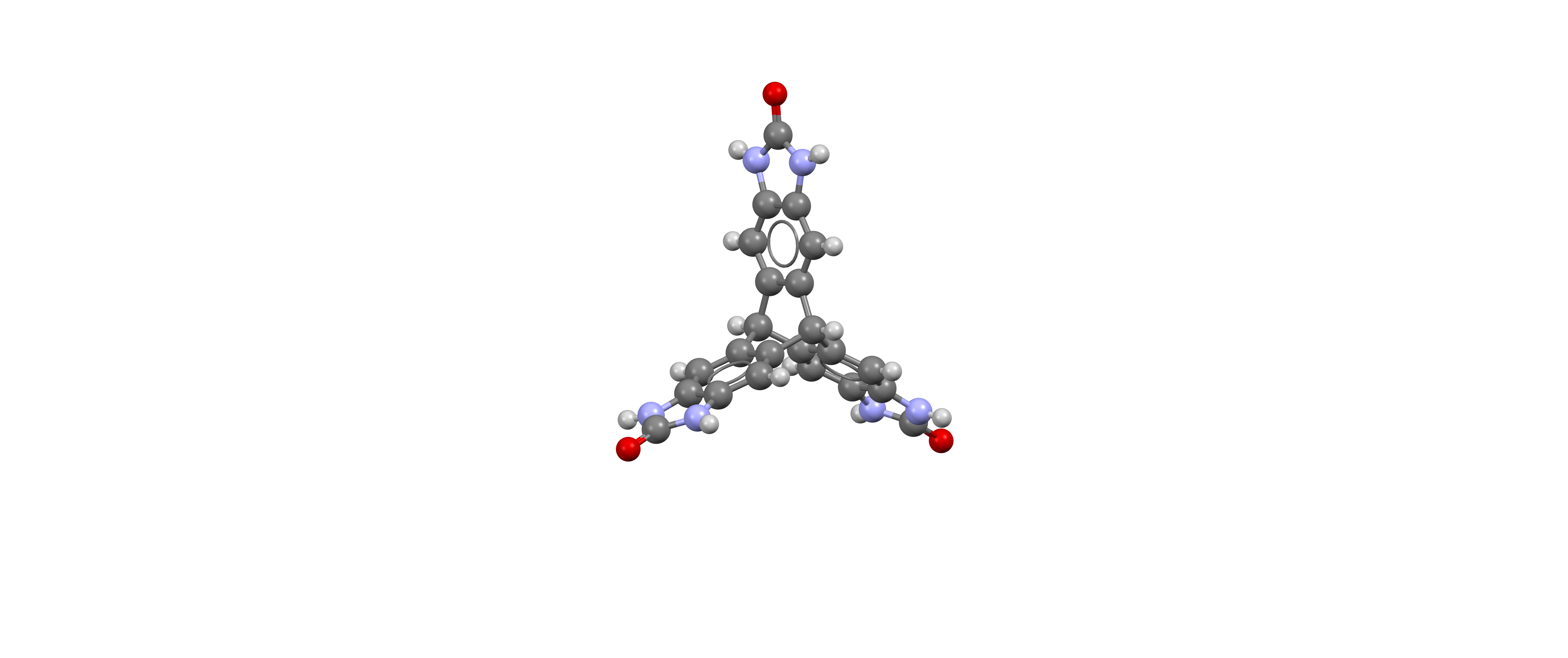}
    \includegraphics[height=32mm]{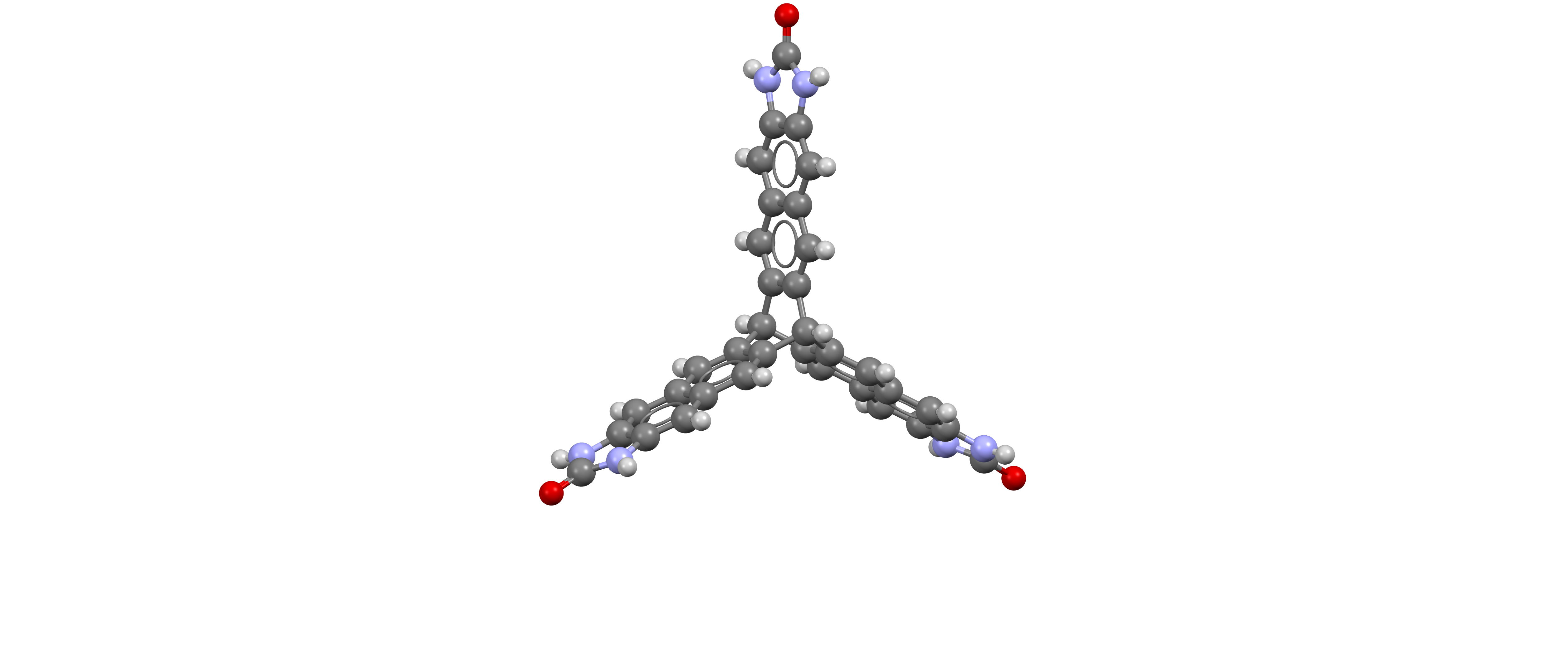}
\caption{T0, T1, T2, and T2E molecules in the four CSP datasets in this section.}    \label{fig:T-molecules}
\end{figure}

In these cases, all experimental crystals have an asymmetric unit consisting of a
single molecule, hence $\CIA=0$ for all versions, which confirms the symmetry
principle saying that real crystals tend to be highly symmetric.
All simulated crystals in the four CSP datasets are based on one of the four
molecules in Fig.~\ref{fig:T-molecules}.

Since each molecule has a rigid shape of three symmetric `arms', its position in
$\R^3$ is uniquely determined by 3 base points at the ends of these `arms' that are
most distant from the molecular centre.
We selected the following 3 base points for each molecule.
T0: mid-points defined by 3 pairs of the most distant carbons from the centre.
T1: three nitrogens.
T2 and T2E: three oxygens.
All values of $\CIA$s in this section were computed on periodic sets obtained by
replacing each molecule with its three base points.
The alternative option of considering all atoms is slower and unnecessary in these
cases, because three base points per molecule suffice to completely reconstruct every
crystal based on one of the molecules T0, T1, T2, and T2E in
Fig.~\ref{fig:T-molecules}.
\smallskip

\begin{figure}
  \centering
  \includegraphics[width=\textwidth]{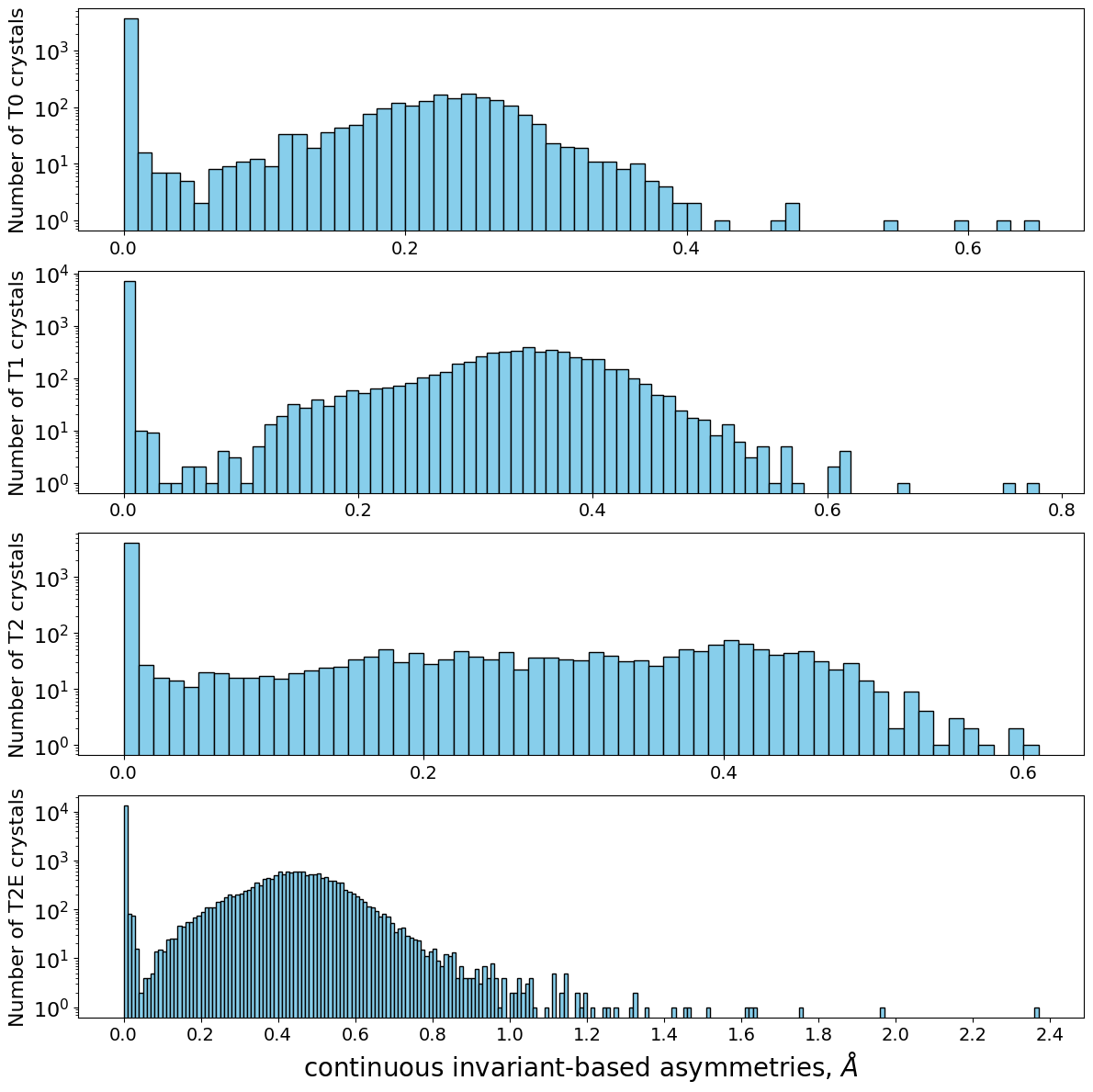}
\caption{The histograms of $\CIA$ for simulated crystals represented by 3 base points
at `ends' of molecules in Fig.~\ref{fig:T-molecules}.
\textbf{Row 1}: T0.
\textbf{Row 2}: T1.
\textbf{Row 3}: T2.
\textbf{Row 4}: T2E. 
}
\label{fig:T_CIA}
\end{figure}

Fig.~\ref{fig:T_CIA} has four histograms of the default $\CIA$ across four CSP
datasets.
In each histogram, the vertical $y$-axis shows the number of crystal structures on
the log scale (as powers of 10) whose $\CIA$s fall in a bin of size $0.01\angstrom$.
The first vertical bin with $\CIA=0$ represents all crystals with $\CIA=0$.
Since any $\CIA$ in Definition~\ref{dfn:CIA} is a min-max or an average of
non-negative distances, all versions of $\CIA$s vanish simultaneously.
\smallskip

All structures in the four CSP datasets were generated with $Z'=1$.
The last stage of energy minimisation allowed this symmetry to be broken, which
explains many cases of $Z'>1$ in Table~\ref{tab:T-crystals_CIA}.
If the generation stage included structures with $Z'\geq 2$, optimised crystals might
have different distribution of $\CIA$s than in Fig.~\ref{fig:T_CIA}.
\smallskip

In appendix~\ref{sec:extra}, Fig.~\ref{fig:T_CIA_inf} contains histograms of
$\CIA_\infty$ based on the $\EMD$ with the ground metric $d_\infty$ in
Definition~\ref{dfn:EMD}.
The Chebyshev metric $d_\infty$ captures the largest deviations, while $d=\RMS$
averages over $k=100$ adjusted inter-atomic distances, $\CIA_\infty$ has a larger
range in comparison with $\CIA$, see maximum values in Table~\ref{tab:T-crystals_CIA}.
\medskip

\begin{table}
\label{tab:T-crystals_CIA}
\centering
\caption{Statistics of $\CIA$ values for the four CSP datasets from
\cite{pulido2017functional}.
The last rows contain Person correlations $r(x,y)$ between energy, density, and new
$\CIA$s.}
\medskip

\begin{tabular}{l|cccc}
CSP datasets & T0 crystals & T1 crystals & T2 crystals & T2E crystals \\ 
\# crystals: all & 5645 & 12524 & 5679 & 29908 \\ \hline
\# crystals: $\CIA\geq 0.001$ & 2024  & 5363 & 1687 & 16491 \\
percentage: $\CIA\geq 0.001$ & 35.8\% & 42.8\% & 29.7\% & 55.1\% \\
maximum $\CIA$, $\angstrom$ & 0.642 & 0.779 & 0.605 & 2.364 \\
\hline
$r(\text{energy}, \text{density})$ & $-0.909$ & $-0.639$ & $-0.377$ & $-0.500$ \\
$r(\text{energy},\CIA)$ & $-0.394$ & $-0.202$ & $+0.022$ & $-0.026$ \\
$r(\text{density},\CIA)$ & $+0.317$ & $+0.148$ & $+0.040$ & $-0.021$
\end{tabular}
\end{table}

\begin{figure}
    \centering
    \includegraphics[width=\textwidth]{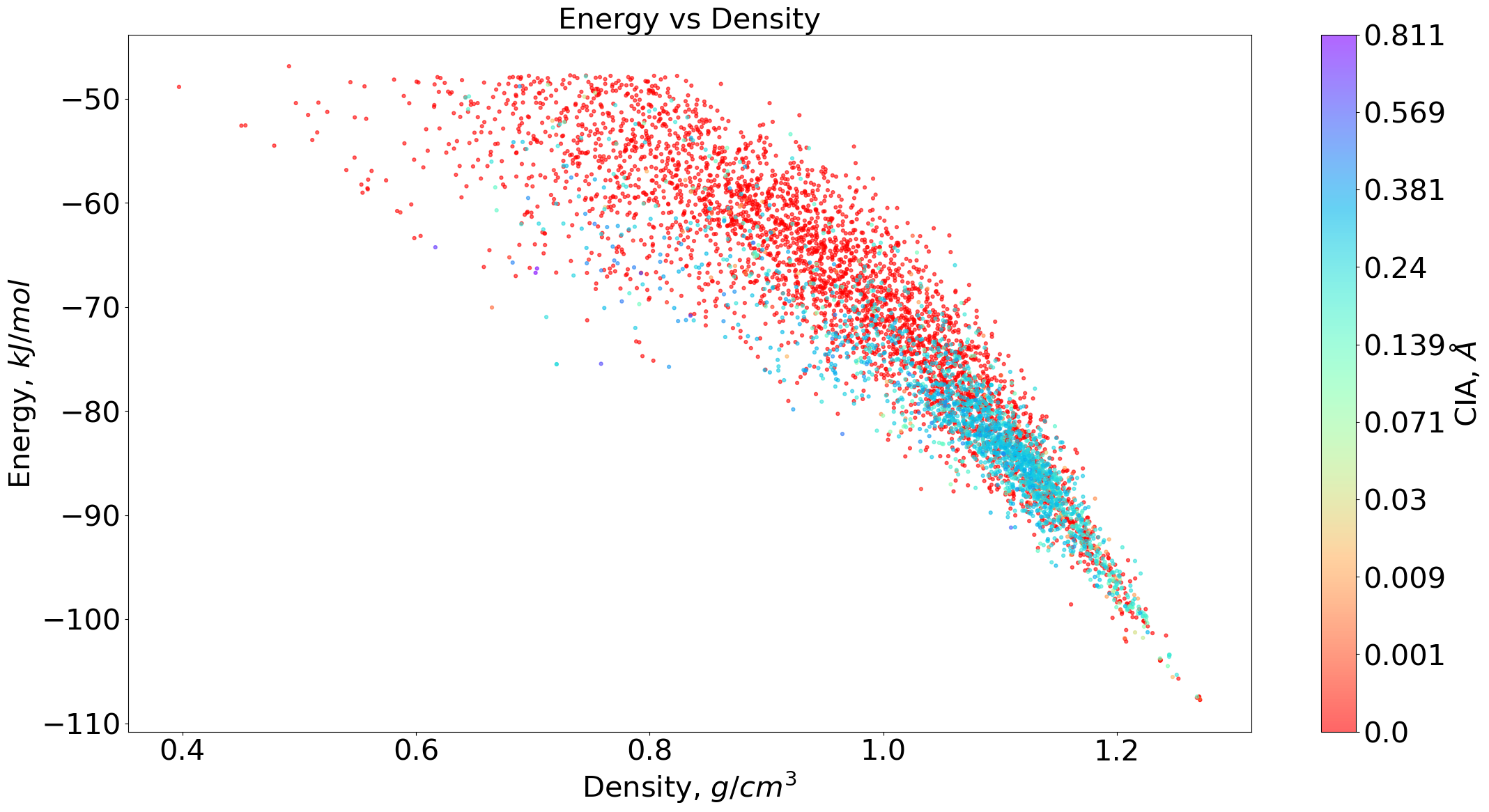}
    \caption{Energy vs density for simulated T0 crystals, coloured by their $\CIA$.}
    \label{fig:T0_energy_vs_density}
\end{figure}

\begin{figure}
    \centering
    \includegraphics[width=\textwidth]{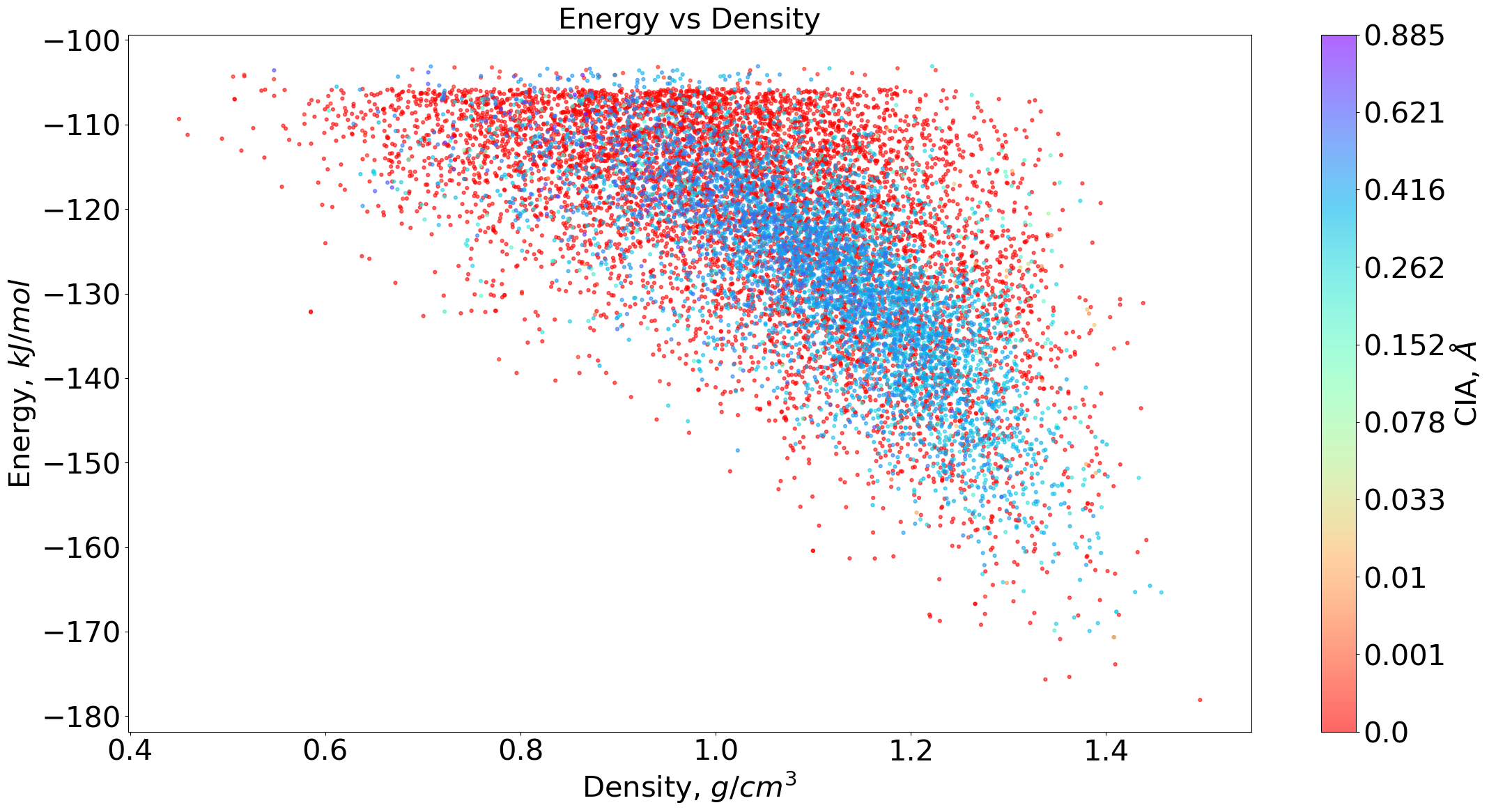}
    \caption{Energy vs density for simulated T1 crystals, coloured by their $\CIA$.}
    \label{fig:T1_energy_vs_density}
\end{figure}

\begin{figure}
    \centering
    \includegraphics[width=\textwidth]{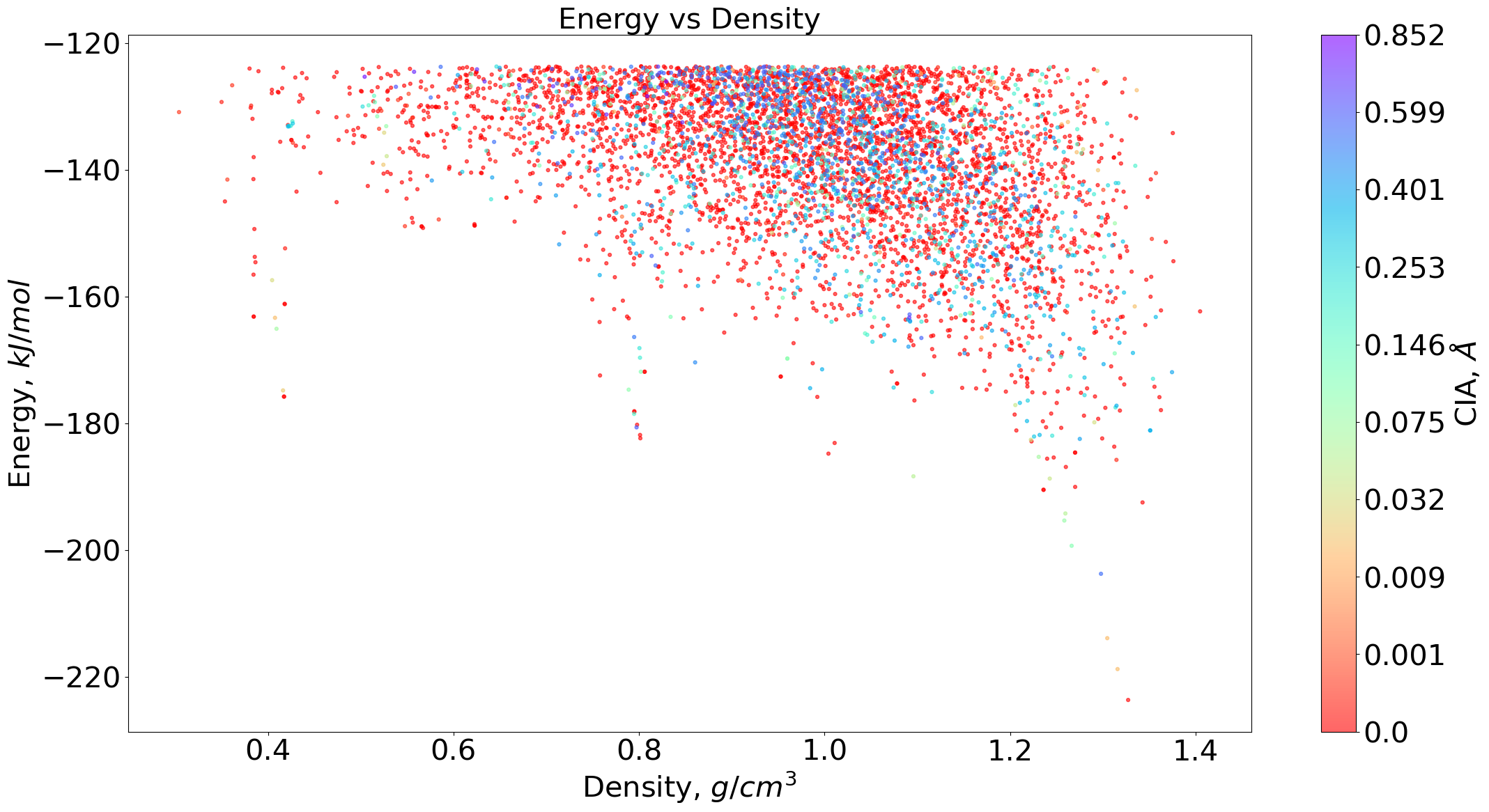}
    \caption{Energy vs density for simulated T2 crystals, coloured by their $\CIA$.}
    \label{fig:T2_energy_vs_density}
\end{figure}

\begin{figure}
    \centering
    \includegraphics[width=\textwidth]{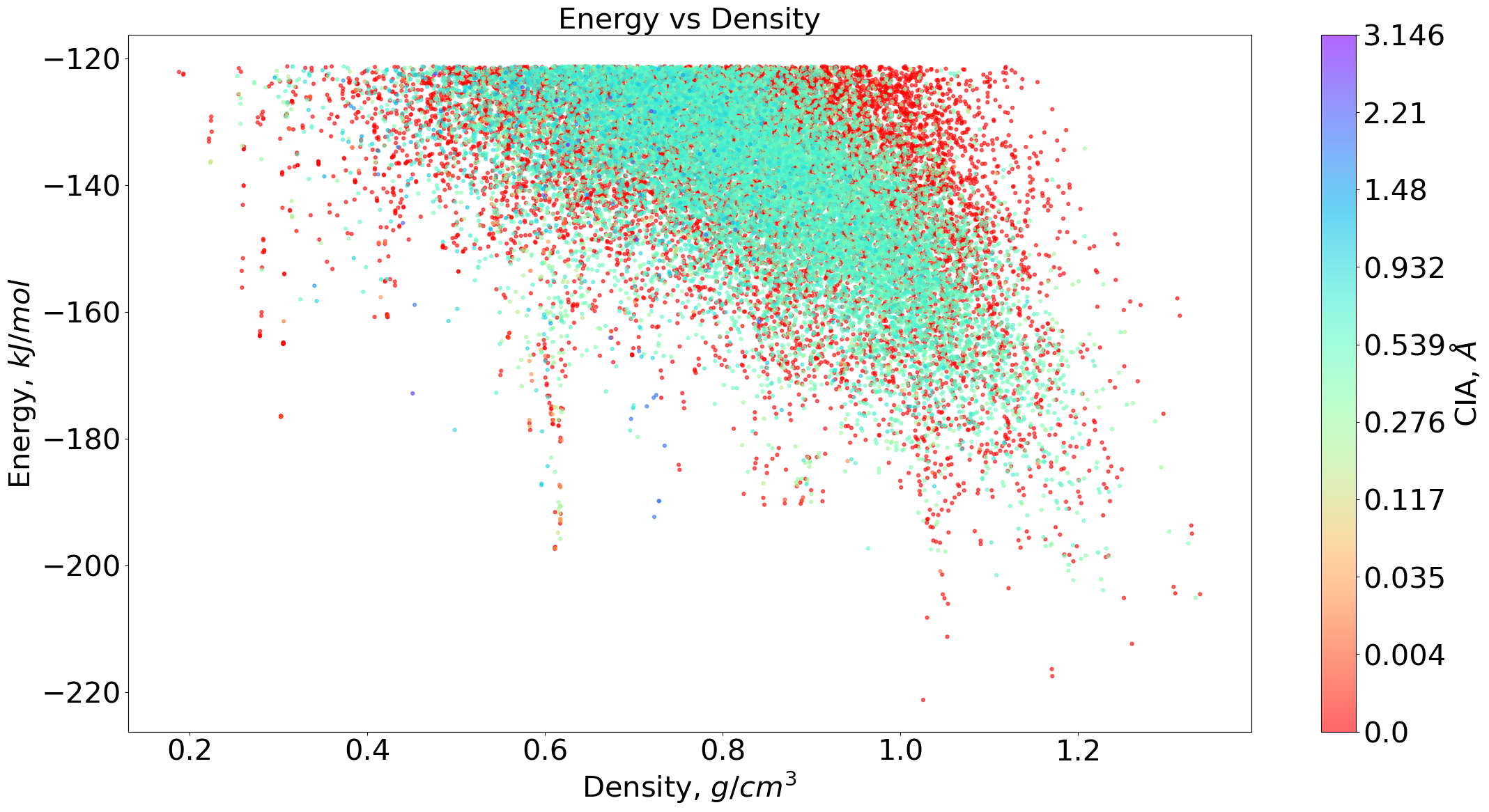}
\caption{Energy vs density for simulated T2E crystals, coloured by their $\CIA$.}    \label{fig:T2E_energy_vs_density}
\end{figure}

CSP datasets are often visualised via energy-density plots, because density is a fast and continuous invariant.
Moreover, density usually indicates stability, because stable crystals tend to be dense.
Figures~\ref{fig:T0_energy_vs_density},~\ref{fig:T1_energy_vs_density},~\ref{fig:T2_energy_vs_density},~\ref{fig:T2E_energy_vs_density}
show these energy-density plots, where each crystal is represented by a point
(density, energy), coloured according to its $\CIA$.
The colour bars on the right-hand side of the plots show the $\CIA$ range, with the bright red colour corresponding to high-symmetry structures with $\CIA=0$.
\smallskip

Table~\ref{tab:T-crystals_CIA} highlights that large subsets (between 30\% and 55\%)
of each CSP dataset have $\CIA>0$.
Since all experimental crystals based on these molecules have $\CIA=0$, all
non-symmetric crystals with $\CIA>0$ are likely non-ideal approximations to symmetric synthesised crystals.
Indeed, if all non-red dots are removed from
Figures~\ref{fig:T0_energy_vs_density},~\ref{fig:T1_energy_vs_density},~\ref{fig:T2_energy_vs_density},
~\ref{fig:T2E_energy_vs_density}, the remaining red dots will still form roughly similar landscapes with all ``minimal spikes'' of density represented by only symmetric crystals with $\CIA=0$ in red.
\smallskip

The Pearson correlation $r(\text{energy},\text{density})$ in Table~\ref{tab:T-crystals_CIA} reflects the inverse dependence on density, because denser crystals tend to be more stable and have lower energies.
This inverse correlation is the strongest with $r=-0.909$ for crystals based on the
smaller molecule T0 and is still noticeable for crystals based on the larger
molecules T1, T2, and T2E.
The new asymmetries $\CIA$ and $\CIA_\infty$ are independent of density and energy due to their low correlations, especially for the T2 and T2E datasets.
All experimental crystals based on these molecules have $\CIA=0$, but their closest
simulated analogues may not have the lowest energies as for the nanoporous T2-$\ga$.

\begin{figure}
    \centering
    \includegraphics[width=\textwidth]{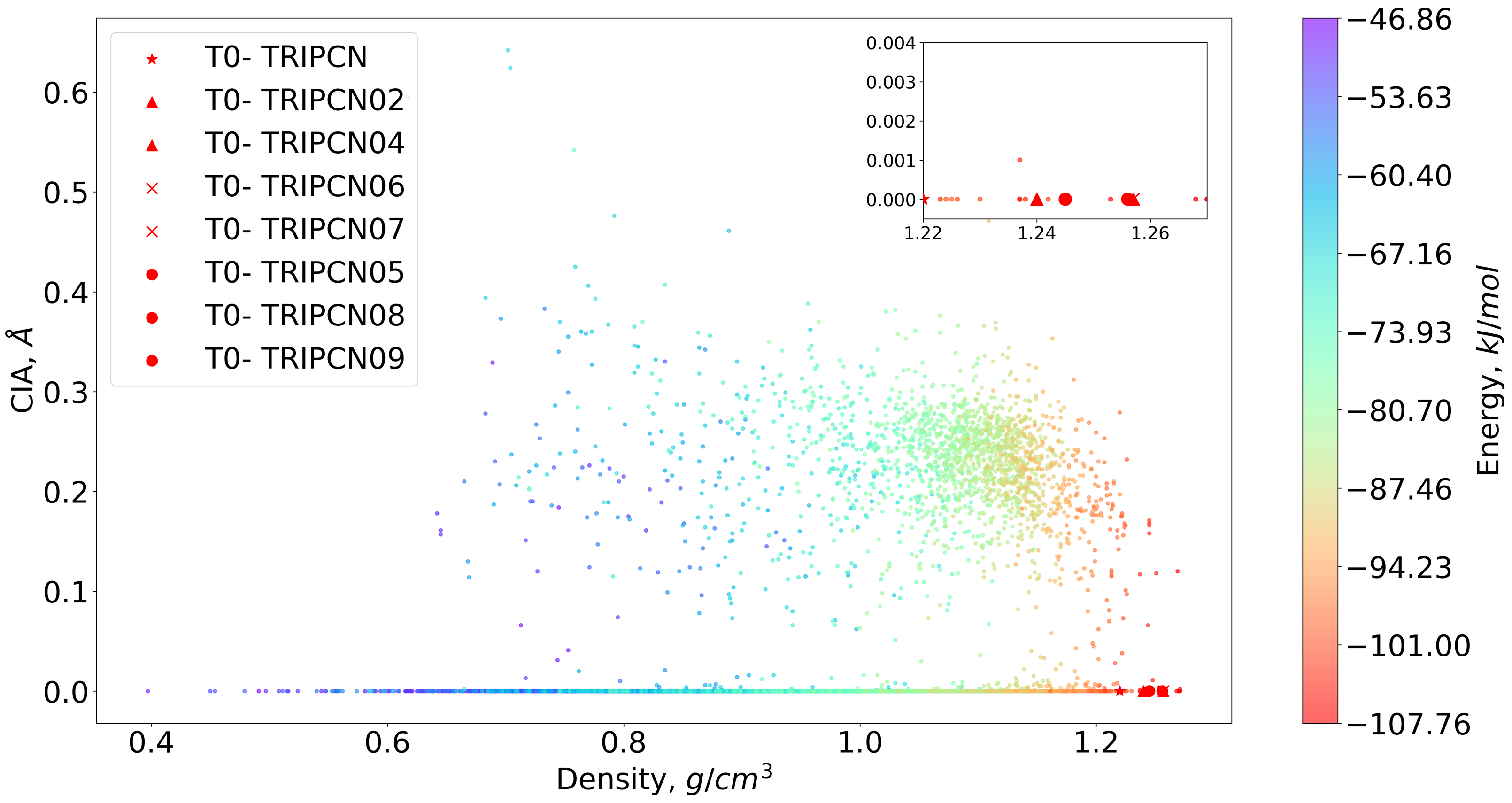}
\caption{$\CIA$ vs density for simulated and experimental T0 crystals in the CSD.}    \label{fig:T0_CIA_vs_density}
\end{figure}

\begin{figure}
    \centering
    \includegraphics[width=\textwidth]{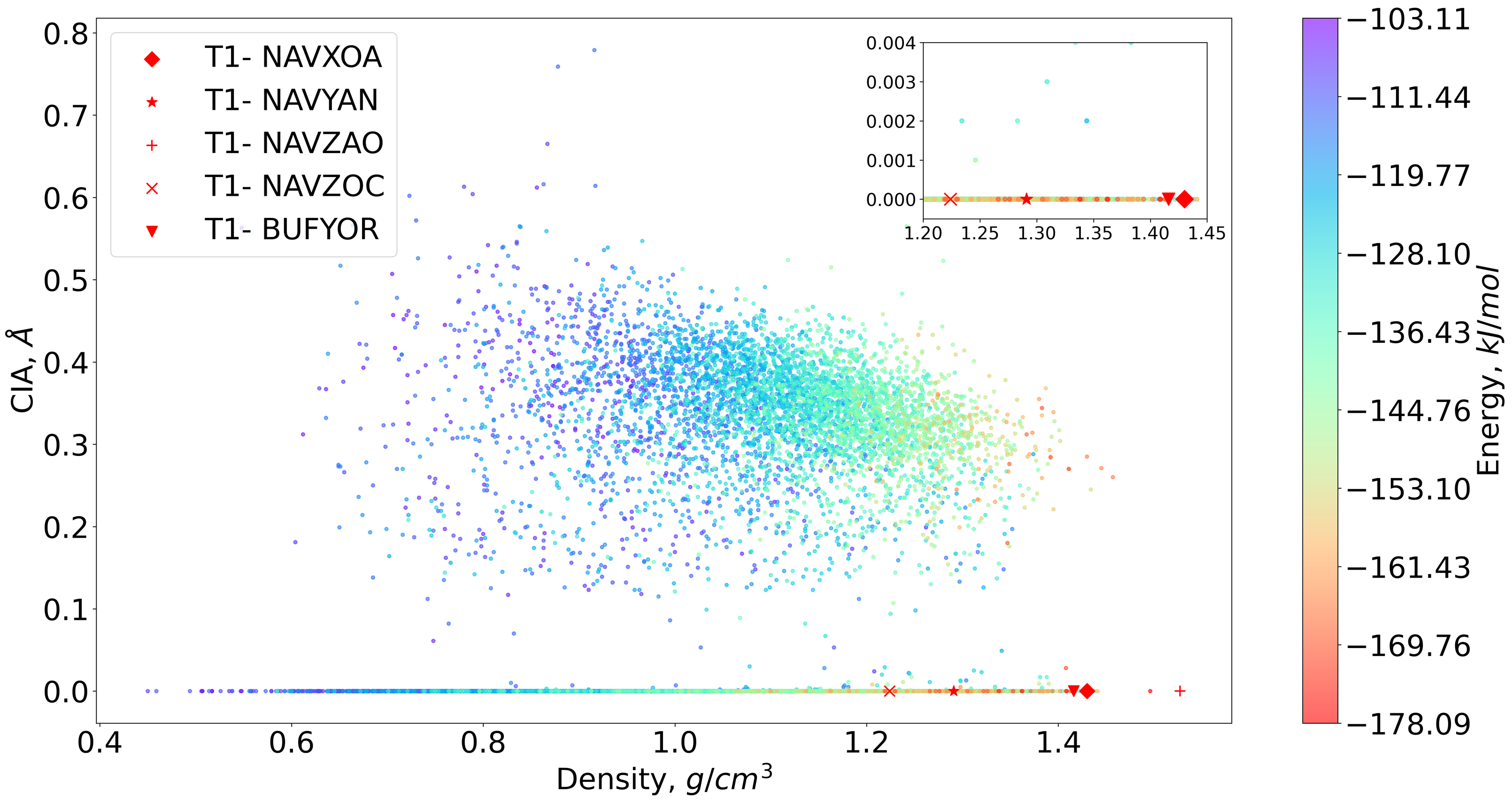}
\caption{$\CIA$ vs density for simulated and experimental T1 crystals in the CSD.}    \label{fig:T1_CIA_vs_density}
\end{figure}

\begin{figure}
    \centering
    \includegraphics[width=\textwidth]{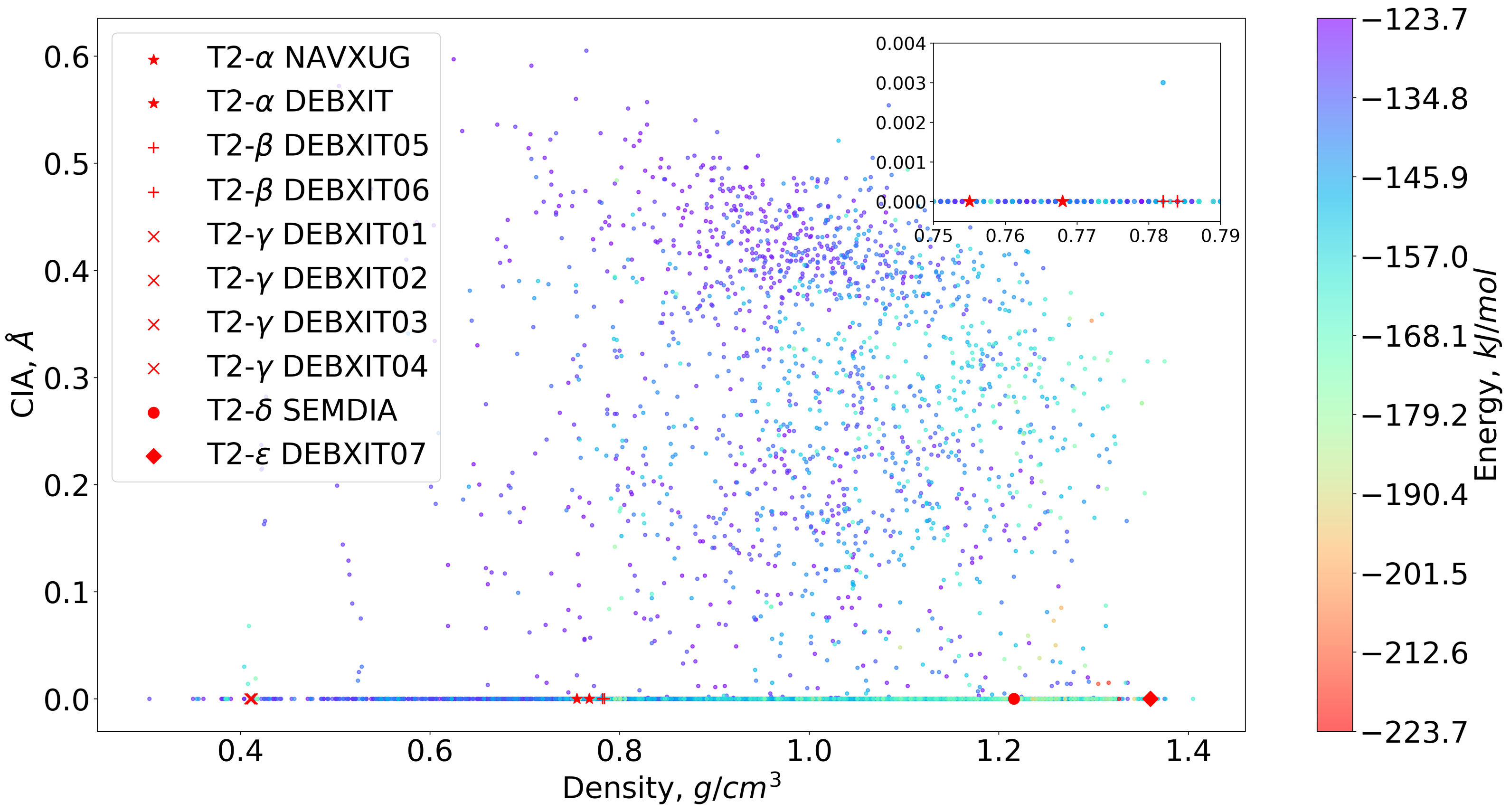}
\caption{$\CIA$ vs density for simulated and experimental T2 crystals in the CSD.}    \label{fig:T2_CIA_vs_density}
\end{figure}

\begin{figure}
    \centering
    \includegraphics[width=\textwidth]{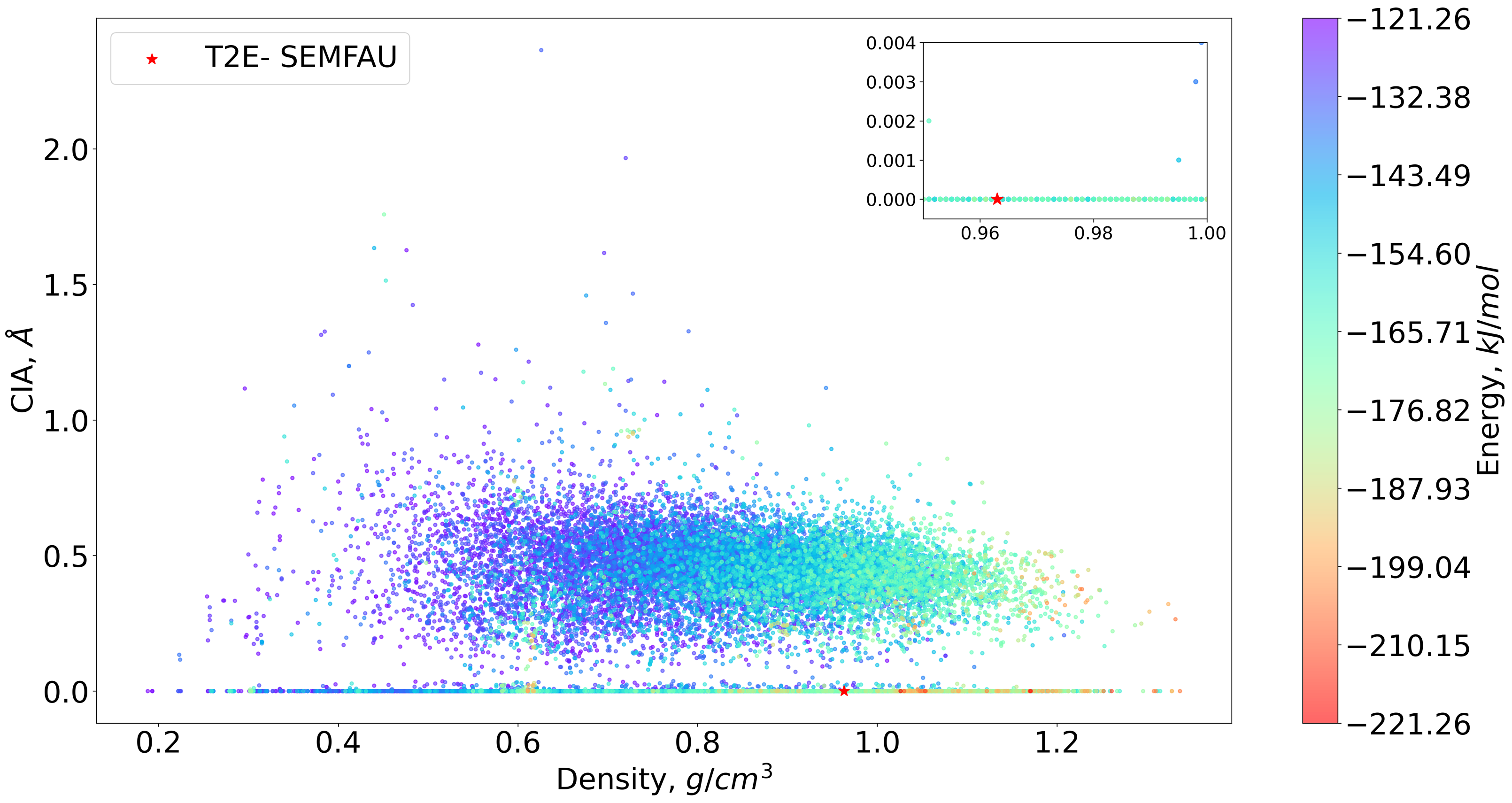}
\caption{$\CIA$ vs density for simulated and experimental T2E crystals in the CSD.}    \label{fig:T2E_CIA_vs_density}
\end{figure}

Figures~\ref{fig:T0_CIA_vs_density},~\ref{fig:T1_CIA_vs_density},~\ref{fig:T2_CIA_vs_density},~\ref{fig:T2E_CIA_vs_density}
show experimental crystals by red marks of various shapes in the coordinates
(density, $\CIA$), and indicate their apparent independence.
In each figure, the top right corner includes a zoomed-in image containing
experimental crystals that are closest by density.
Though many simulated crystals are symmetric with $\CIA=0$, all non-symmetric
crystals form noisy clouds with variable energies.
The visible gaps between these clouds and the horizontal axis $\CIA=0$ confirm a
local version of the symmetry principle saying that a nearly symmetric structure
likely converges to a higher symmetry version with $\CIA=0$.
\smallskip

\begin{figure}
\centering
\includegraphics[width=\textwidth]{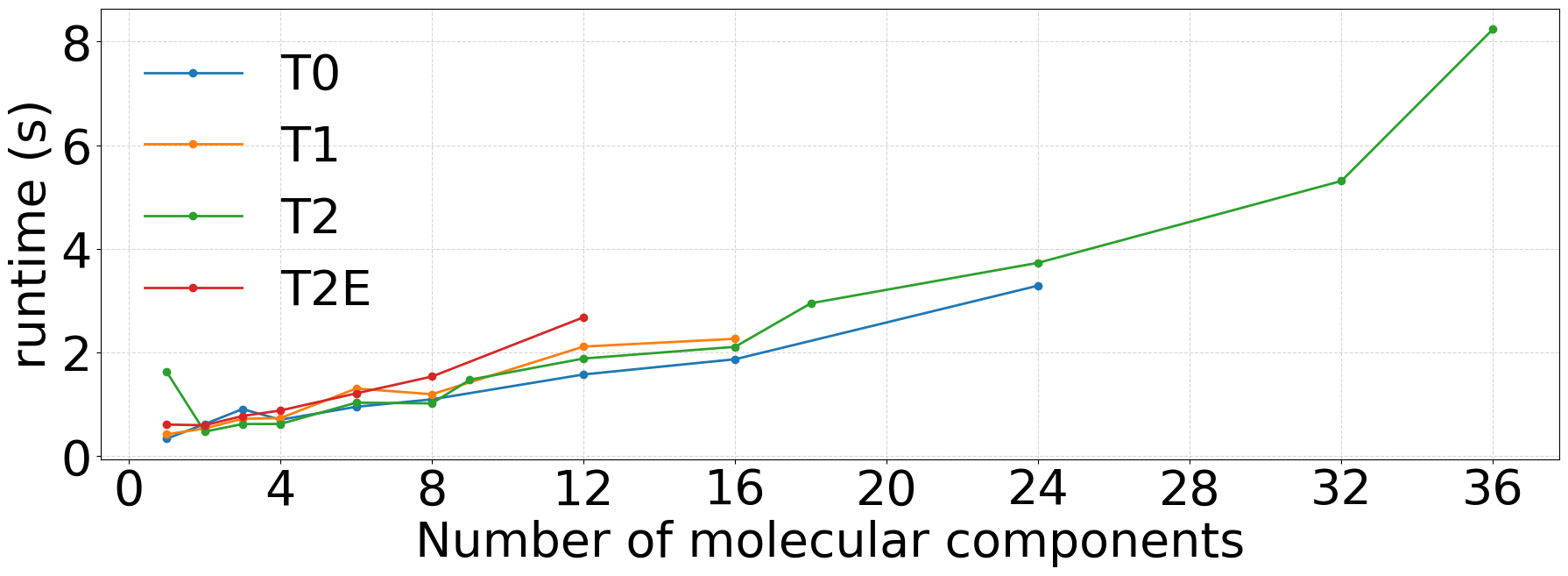}
\caption{Average running times (in seconds) of $\CIA$ on four CSP datasets vs the
number $G$ of molecules in asymmetric units, performed on a modest machine with CPU
13th Gen Intel(R) Core(TM) i7-1355U (1.70 GHz) and RAM 16GB.
}
\label{fig:CIA_T-crystals_times}
\end{figure}

Figure~\ref{fig:CIA_T-crystals_times} shows the average running times vs the number
$Z$ of molecular components in a unit cell.
This number $Z$ goes up to 36 and coincides with $G$, because all finally optimised
crystals are saved in the simplest translation group P1.

\section{Continuous asymmetries of all experimental crystals in the CSD}
\label{sec:experimental}

This section describes a large-scale analysis of asymmetries in the whole CSD. 
Each crystal is represented by a periodic set of all its atoms.
We considered all periodic crystals with complete 3D geometry, no disorder, and based
on a chemically unique molecule.
Though Definition~\ref{dfn:CIA} can be applied to geometric blocks of different
sizes, we postpone the more complicated case of co-crystals to future work.
\smallskip

\begin{figure}
\centering
\includegraphics[width=\textwidth]{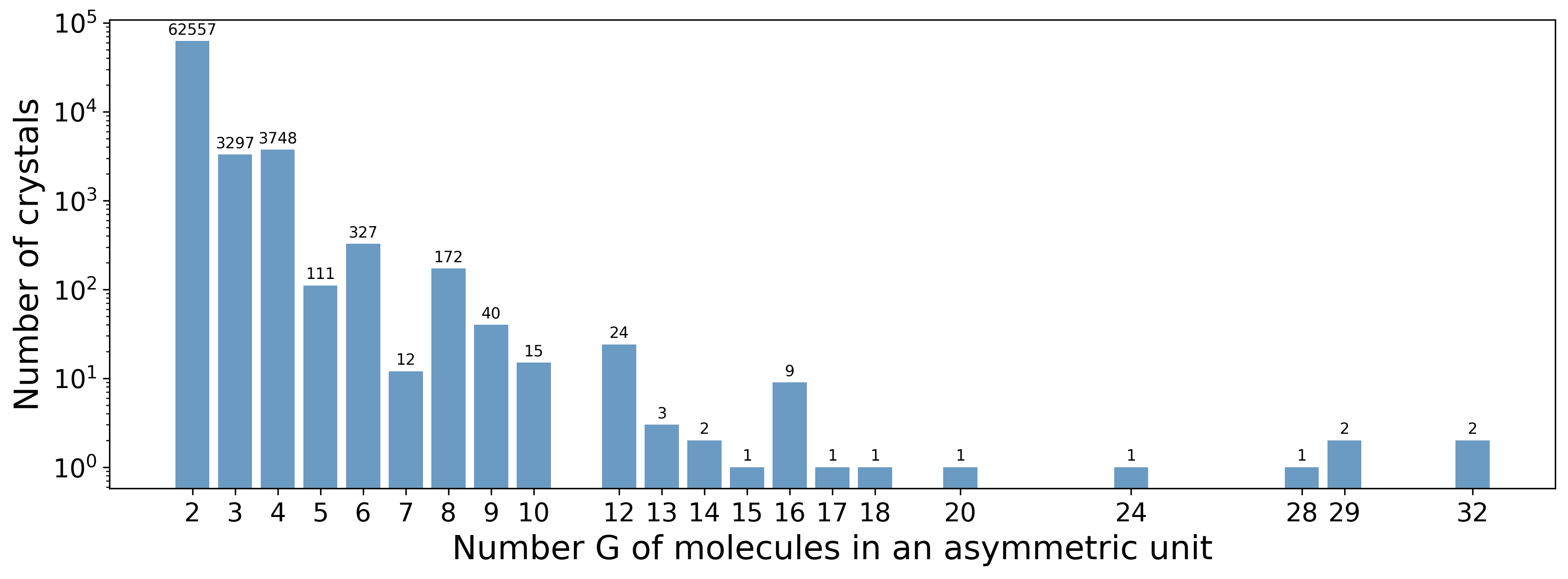}
\caption{The histogram of integer numbers $G$ for all 69,196 periodic crystals in the
CSD that have $G\geq 2$ chemically equivalent blocks in their asymmetric units.}
\label{fig:CSD_Gblocks}
\end{figure}

\begin{figure}
\centering
\includegraphics[width=\textwidth]{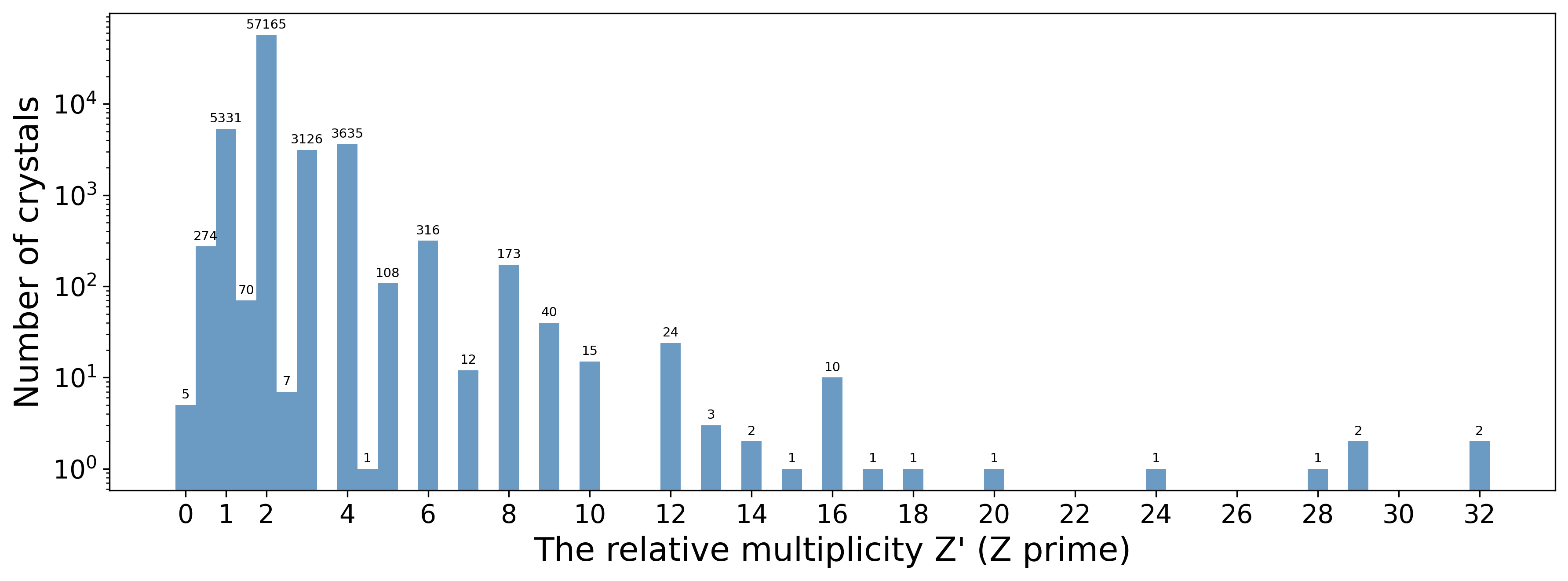}
\caption{The histogram of 
$Z'$ with bin size 0.5 for all 69,196 periodic crystals in the CSD that have $G\geq
2$ chemically equivalent blocks in their asymmetric units.}
\label{fig:CSD_Zprime}
\end{figure}

The snapshot of the CSD on 12th November 2025 contained 1,394,755 entries, including
907223 crystals without disorder.
Among them, 69,196 crystals have asymmetric units containing $G\geq 2$ molecules that
all have the same composition, where $G$ was computed by the CSD Python API as the
number of components in the list crystal.asymmetric\_unit\_molecules.
Some crystals with the highest $Z'$ values from https://zprime.co.uk/database, such
as OGUROZ ($Z'=56$), TMESNH ($Z'=32$), IDOSID ($Z'=24$), and VIFXEQ ($Z'=24$), were
excluded because of disorder.
\smallskip

\begin{figure}
\centering
\includegraphics[width=\textwidth]{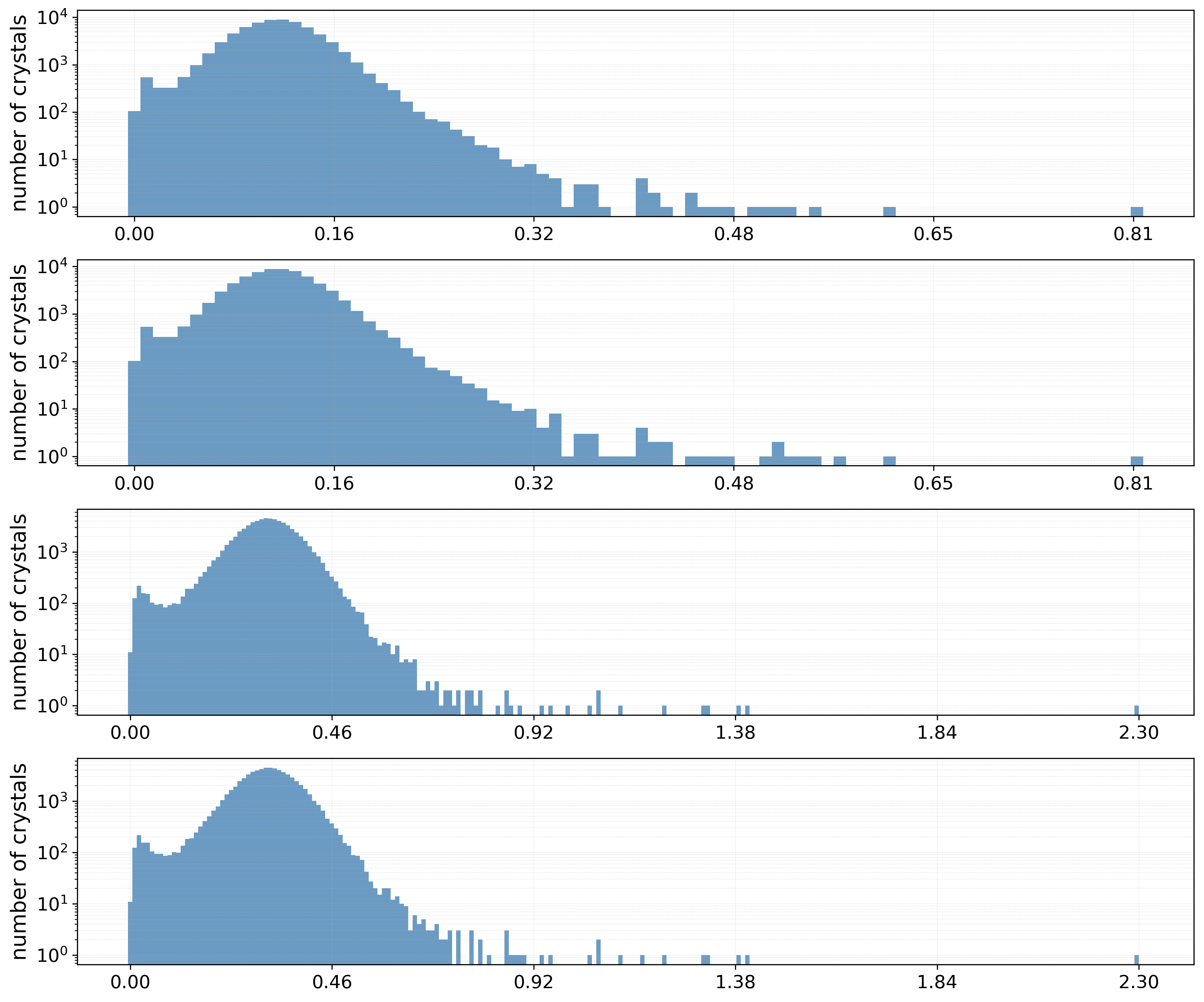}
\caption{The histograms of $\CIA$s on the log scale with bin size $0.01\angstrom$ for
all 69,196 periodic crystals in the CSD that have $G\geq 2$ chemically equivalent
molecules in asymmetric units.
\textbf{Row 1}: $\CIA$.
\textbf{Row 2}: $\aCIA$.
\textbf{Row 3}: $\CIA_\infty$.
\textbf{Row 4}: $\aCIA_\infty$.
}
\label{fig:CSD_CIAs}
\end{figure}

Figures~\ref{fig:CSD_Gblocks}, \ref{fig:CSD_Zprime}, and \ref{fig:CSD_CIAs} show the
histograms of $G,Z'$, and four $\CIA$s for this subset of the CSD, respectively,
where $Z'$ was computed as crystal.z\_prime by the CSD Python API.
The number $Z$[CIF] of molecules in the full motif was taken from lines
``\_cell\_formula\_units\_Z'' in CIFs from the CSD, which sometimes differs from
$Z$[CSD], computed as the number of components in the list crystal.molecule. 
CSD Python API.
\smallskip

Table~\ref{tab:CSD_Z'_extreme} shows all four versions of $\CIA$s for the most
extreme crystals in the CSD: five crystals with the lowest $Z'\leq 0.33$ and five
crystals with the highest $Z'\geq 28$.

\begin{table}
\centering
\caption{$\CIA$s of the crystals with the lowest and largest relative multiplicities
in the CSD.
The numbers $Z$ and $G$ count molecules in a unit cell and an asymmetric unit,
respectively.
}
\label{tab:CSD_Z'_extreme}
\medskip

\begin{tabular}{lcccccccc}
CSD id & $Z$[CIF] & $G$ blocks & $Z'$[CSD] & $\CIA, \angstrom$ & $\aCIA, \angstrom$ &
$\CIA_\infty, \angstrom$ & $\aCIA_\infty, \angstrom$ \\
\hline
VESWEZ & 2 & 2 & 0.083 & 0.204 & 0.204 & 0.580 & 0.580 \\
ELIQIZ02 & 3 & 2 & 0.083 & 0.226 & 0.226 & 0.807 & 0.807 \\
ZOKYEH01 & 16 & 2 & 0.167 & 0.086 & 0.086 & 0.241 & 0.241 \\
RARTEK & 16 & 2 & 0.17 & 0.125 & 0.125 & 0.332 & 0.332 \\
ZAVJOV & 2 & 2 & 0.33 & 0.090 & 0.090 & 0.241 & 0.241 \\
\hline
QILJII01  & 112 & 28 & 28 & 0.168 & 0.185 & 0.397 & 0.426 \\
LOFRAD & 116 & 29 & 29 & 0.149 & 0.183 & 0.420 & 0.499 \\
LOFRAD01 & 116 & 29 & 29 & 0.149 & 0.185 & 0.434 & 0.506 \\
JIPTIL09 & 32 & 32 & 32 & 0.104 & 0.109 & 0.266 & 0.282 \\
JIPTIL10 & 32 & 32 & 32 & 0.102 & 0.109 & 0.265 & 0.282 \\
\end{tabular}
\end{table}

\begin{figure}
\centering
\includegraphics[width=\textwidth]{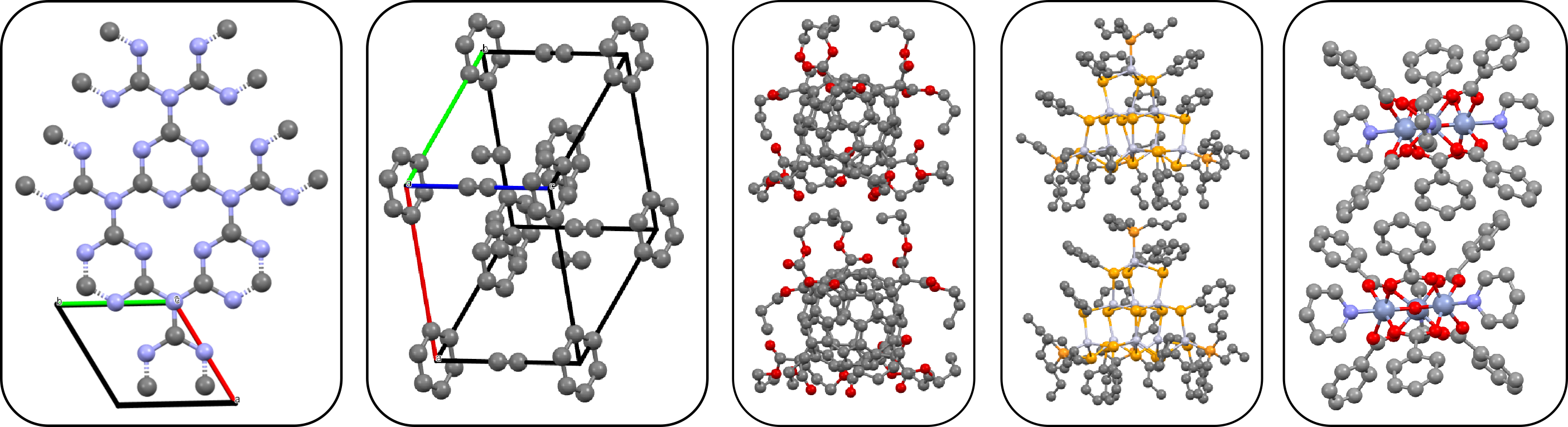}
\caption{The crystals with the lowest $Z'$ from Table~\ref{tab:CSD_Z'_extreme} shown
without hydrogen atoms.
\textbf{1st}: VESWEZ. 
\textbf{2nd}: ELIQIZ02. 
\textbf{3rd}: ZOKYEH01. 
\textbf{4th}: RARTEK. 
\textbf{5th}: ZAVJOV. 
}
\label{fig:CSD_Z'_lowest}
\end{figure}

In Table~\ref{tab:CSD_Z'_extreme}, crystal VESWEZ has $G=2$ components \ce{CN2} in
geometrically non-equivalent positions: in one \ce{CN2}, both nitrogens are linked to
two carbons; in another \ce{CN2}, the two nitrogens are linked to 2 and 3 carbons,
see Fig.~\ref{fig:CSD_Z'_lowest}.
Crystal ELIQIZ02 has molecules \ce{C6H6} and \ce{C2H2}, and its asymmetric unit
consists of $G=2$ geometrically different carbons: one from \ce{C6H6} and another
from \ce{C2H2}.
Crystal ZOKYEH01 consists of a big molecule of \ce{C60} with extra tails, but its
asymmetric unit was also split into $G=2$ blocks \ce{C10O2}, which apparently have
non-isometric positions within the full crystal.
Crystal RARTEK and ZAVJOV similarly consist of big molecules based on $G=2$ blocks in
asymmetric units, whose positions can not be matched by isometry preserving the whole
crystal.
The last three cases show that molecular crystals will benefit from quantifying
asymmetry at the level of full molecules, because asymmetric units may not split into
uniquely defined molecules or geometric blocks of atoms.

\medskip

\begin{table}
\centering
\caption{
CIAs of the well-known polymorphs of
artemisinin (QNGHSU01), pyridine (PYRDNA04), para-chlorophenol ($\al$-form CLPHOL12
and $\be$-form CLPHOL13), and the famous ROY molecule (R05 polymorph QAXMEH31 and R18
polymorph QAXMEH57). }
\label{tab:CSD_famous_polymorphs}

\begin{tabular}{lcccccccc}
CSD id & $Z$[CIF] & $G$ blocks & $Z'$[CSD] & $\CIA, \angstrom$ & $\aCIA, \angstrom$ &
$\CIA_\infty, \angstrom$ & $\aCIA_\infty, \angstrom$ \\
\hline
QNGHSU01 & 4 & 4 & 4 & 0.357 & 0.379 & 1.093 & 1.096 \\
QAXMEH31 & 2 & 2 & 2 & 0.440 & 0.440 & 1.098 & 1.098 \\
QAXMEH57 & 2 & 2 & 2 & 0.807 & 0.807 & 1.602 & 1.602 \\
CLPHOL12 & 2 & 2 & 2 & 0.790 & 0.790 & 2.594 & 2.594 \\
CLPHOL13 & 2 & 2 & 2 & 0.575 & 0.575 & 1.132 & 1.132 \\
PYRDNA04 & 4 & 4 & 4 & 1.971 & 2.096 & 2.756 & 2.861 \\
\end{tabular}
\end{table}

\begin{figure}
\centering
\includegraphics[width=\textwidth]{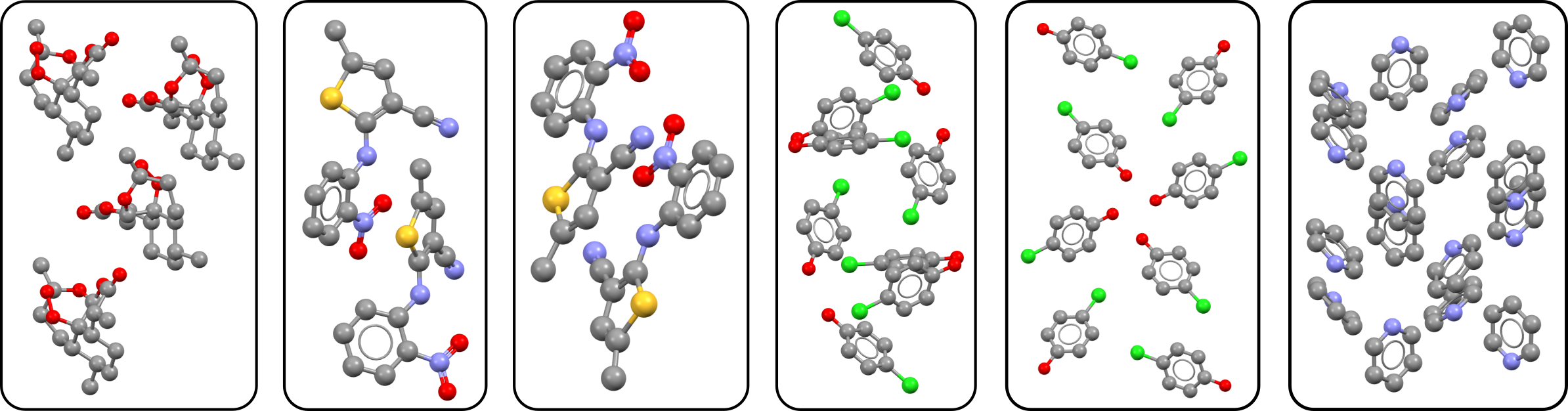}
\caption{Six famous polymorphs whose CIAs are listed in
Table~\ref{tab:CSD_famous_polymorphs}.
From left to right: QNGHSU01, QAXMEH31, QAXMEH57, CLPHOL12, CLPHOL13, PYRDNA04.}
\label{fig:CSD_famous_polymorphs}
\end{figure}

Table~\ref{tab:CSD_CIA_lowest} lists the 10 crystals from with the lowest values of
$\CIA$s.
The first three crystals have $\CIA=0$ with 3 decimal places, so their $Z'\geq 2$
might be corrected.

\begin{table}
\centering
\caption{Ten crystals with the lowest $\CIA$ among 69,196 periodic crystals in the
CSD that have $Z'\geq G\geq 2$ chemically equivalent blocks in their asymmetric
units.
}
\label{tab:CSD_CIA_lowest}
\medskip

\begin{tabular}{lcccccccc}
CSD id & $Z$[CIF] & $G$ blocks & $Z'$[CSD] & $\CIA, \angstrom$ & $\aCIA, \angstrom$ &
$\CIA_\infty, \angstrom$ & $\aCIA_\infty, \angstrom$ \\
\hline
IYIWIY & 8 & 8 & 8 & 0.000 & 0.000 & 0.000 & 0.000 \\
GLYCIN81 & 2 & 2 & 2 & 0.000 & 0.000 & 0.000 & 0.000 \\
YOSNEZ05 & 2 & 2 & 2 & 0.000 & 0.000 & 0.000 & 0.000 \\
\hline
GIBVOG & 2 & 2 & 2 & 0.000 & 0.000 & 0.001 & 0.001 \\
GLYCIN82 & 3 & 3 & 3 & 0.001 & 0.001 & 0.002 & 0.002 \\
KAVXUE & 1 & 2 & 2 & 0.002 & 0.002 & 0.005 & 0.005 \\
ADUWED & 64 & 2 & 2 & 0.002 & 0.002 & 0.006 & 0.006 \\
CINMAC13 & 2 & 2 & 2 & 0.002 & 0.002 & 0.010 & 0.010 \\
XOTRAB & 4 & 2 & 2 & 0.003 & 0.003 & 0.007 & 0.007 \\
COTZES & 6 & 2 & 2 & 0.003 & 0.003 & 0.009 & 0.009
\end{tabular}
\end{table}

The value $\CIA=0$ means that all molecules are geometrically equivalent, i.e. can be
exactly matched by isometry that preserves the whole crystal.
In this case, an asymmetric unit should contain only one molecule ($G=1$), so $Z'\leq
1$ is expected.
\smallskip

\begin{figure}
\centering
\includegraphics[width=\textwidth]{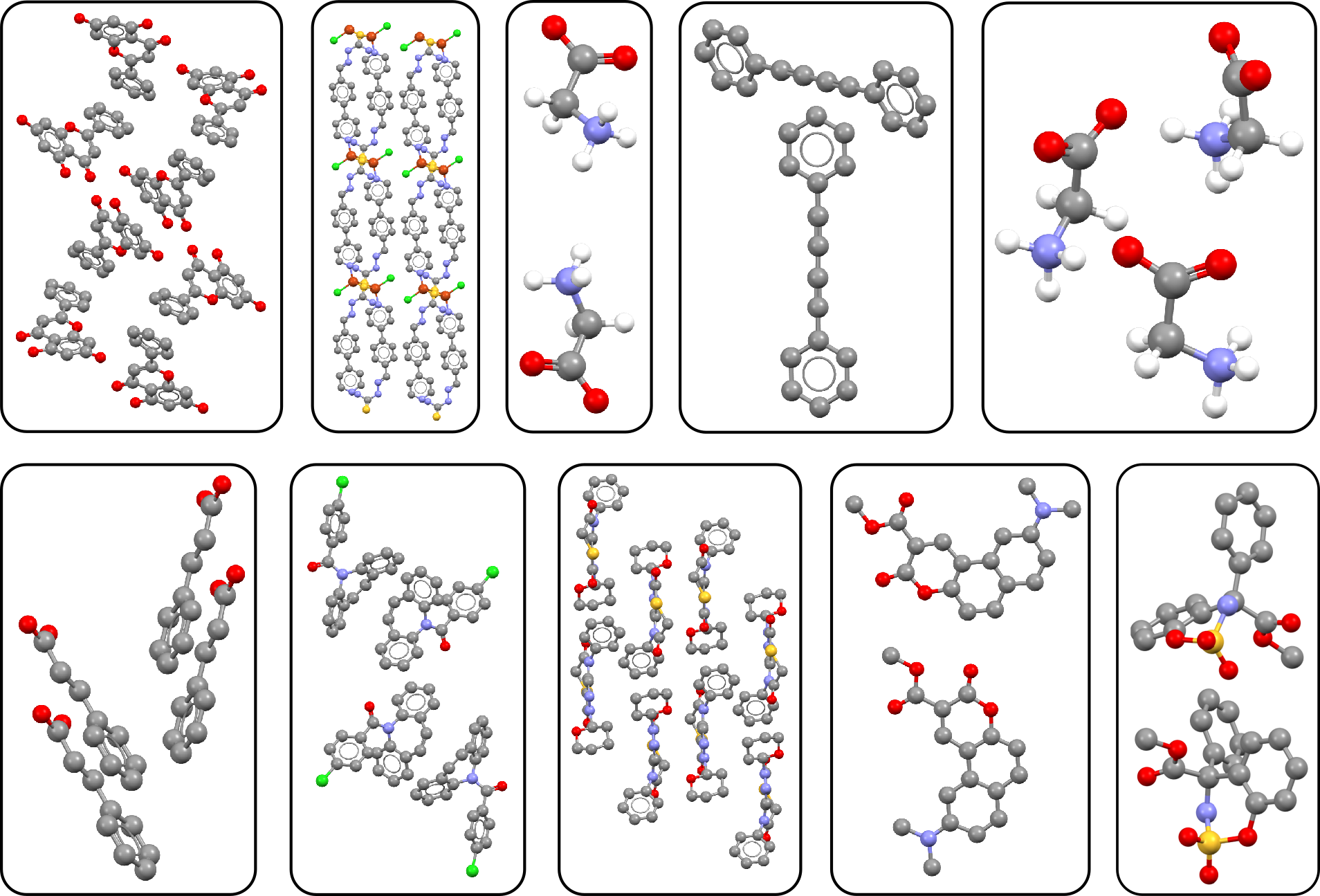}
\caption{Ten crystals (some shown without hydrogen atoms) from
Table~\ref{tab:CSD_CIA_lowest} with very low $\CIA\geq 0$.
\textbf{Top} from left to right: IYIWIY, GIBVOG, GLYCIN81, YOSNEZ05, GLYCIN82. 
\textbf{Bottom}: CINMAC13, KAVXUE, ADUWED, XOTRAB, COTZES.}
\label{fig:CSD_lowCIA}
\end{figure}

One explanation is a potentially wrong space group \cite{henling2014some}.
For example, IYIWIY has the space group P1, but looks more symmetric in the first
picture of Fig.~\ref{fig:CSD_lowCIA}.
Both structures IYIWIY and YOSNEZ05 were obtained from powder data, so their space
groups might need re-checking.
Since all $\CIA$s continuously change under atomic perturbations by
Theorem~\ref{thm:CIA}, there is no need to search for a higher symmetry group, which
drops to the simplest group P1 under almost any noise anyway.
The values of $Z$[CIF] can be corrected for all entries with $Z<G$ in
Tables~\ref{tab:CSD_CIA_lowest} and~\ref{tab:CSD_highZ'_lowCIA}, because a unit cell
should not have fewer molecules than in an asymmetric unit.

\begin{table}
\centering
\caption{Almost symmetric crystals with high values $Z'\geq 5$ but low $\CIA\leq
0.021\angstrom$.}
\label{tab:CSD_highZ'_lowCIA}
\medskip

\begin{tabular}{lcccccccc}
CSD id & $Z$[CIF] & $G$ blocks & $Z'$[CSD] & $\CIA, \angstrom$ & $\aCIA, \angstrom$ &
$\CIA_\infty, \angstrom$ & $\aCIA_\infty, \angstrom$ \\
\hline
TEGBEP & 1 & 6 & 6 & 0.010 & 0.011 & 0.030 & 0.032 \\
HOGKAR & 12 & 6 & 6 & 0.010 & 0.011 & 0.032 & 0.034 \\
GINHIX & 6 & 6 & 6 & 0.011 & 0.012 & 0.034 & 0.039 \\
EVIWUE & 12 & 6 & 6 & 0.012 & 0.014 & 0.051 & 0.059 \\
LEMWOR & 2 & 6 & 6 & 0.013 & 0.013 & 0.040 & 0.043 \\
YIVHER & 10 & 5 & 5 & 0.015 & 0.016 & 0.053 & 0.058 \\
IFOFAN & 10 & 5 & 5 & 0.020 & 0.020 & 0.070 & 0.077 \\
EDUCAL & 12 & 6 & 6 & 0.020 & 0.022 & 0.062 & 0.071 \\
ROTSAY & 18 & 9 & 9 & 0.021 & 0.023 & 0.060 & 0.067 \\
CIDHAB & 1 & 12 & 12 & 0.021 & 0.023 & 0.067 & 0.071 \\
\end{tabular}
\end{table}

In conclusion, the relative multiplicity $Z'$ discontinuously changes under almost
any perturbation, the proposed $\CIA$ in Definition~\ref{dfn:CIA} is continuous by
Theorem~\ref{thm:CIA}.
For the CSP datasets in section~\ref{sec:simulated}, about a half of all 50K+
simulated crystals have $\CIA>0$, while all experimental crystals have $\CIA=0$, see
Table~\ref{tab:T-crystals_CIA}.
Moreover, these continuous and fast asymmetries are not correlated with density and
energy.
The large-scale experiments on the CSD show that many non-symmetric crystals with
high $Z'$ have low $\CIA$s in Table~\ref{tab:CSD_highZ'_lowCIA} and hence are
geometrically close to more symmetric forms.
This work was supported by the Royal Society APEX fellowship "New geometric methods for mapping the space of periodic crystals" (APX/R1/231152) of the last author.

\referencelist[continuous-asymmetry]

\appendix

\section{Extra experimental results for simulated crystals}
\label{sec:extra}

This appendix contains extra plots for other versions of $\CIA$s.  
 
\begin{figure}
  \centering
  \includegraphics[width=\textwidth]{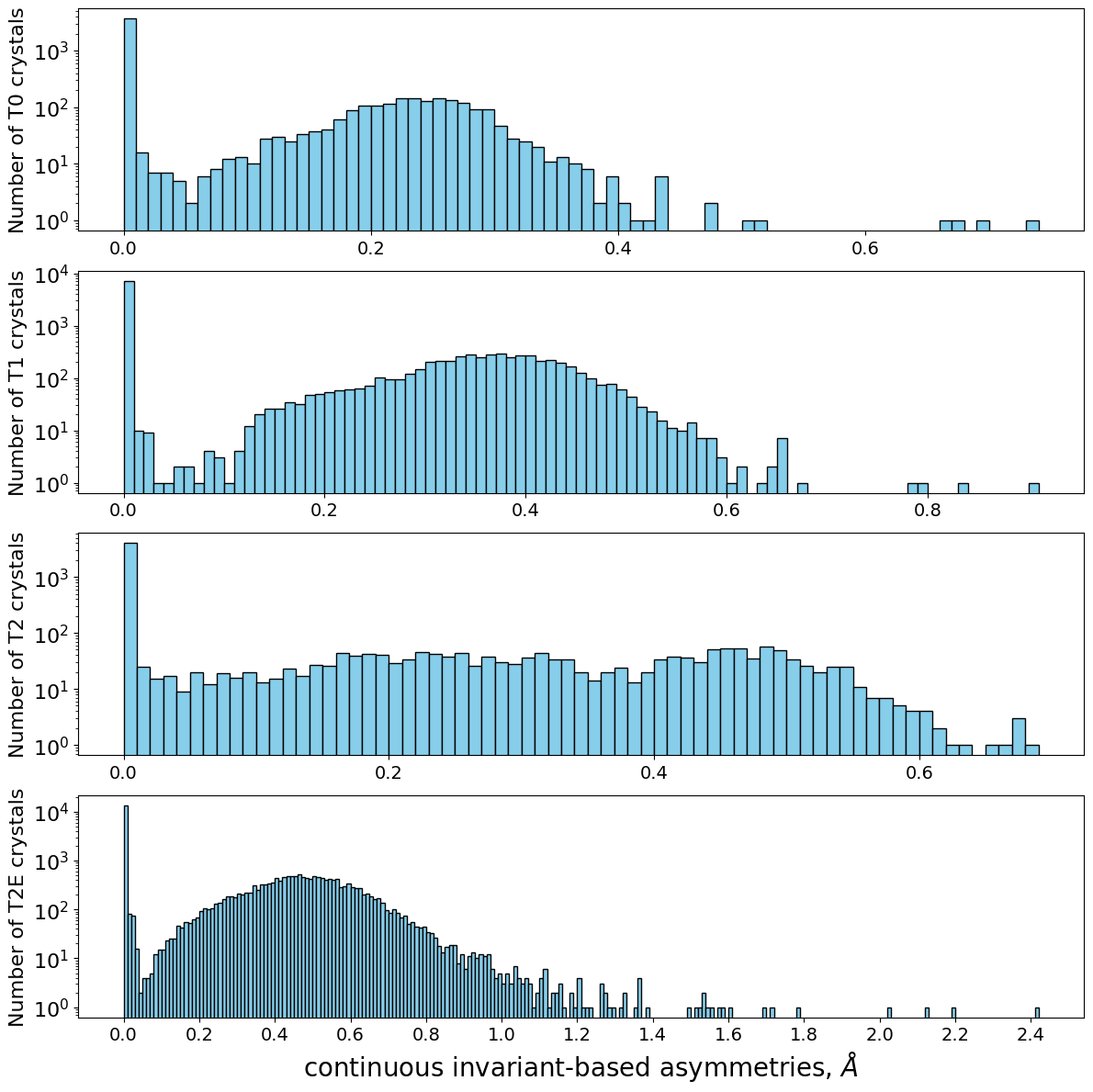}
\caption{The histograms of $\aCIA$ for simulated crystals represented by 3 base
points at `ends' of molecules in Fig.~\ref{fig:T-molecules}.
\textbf{Row 1}: T0.
\textbf{Row 2}: T1.
\textbf{Row 3}: T2.
\textbf{Row 4}: T2E. 
}
\label{fig:T_aCIA}
\end{figure}

\begin{figure}
  \centering
  \includegraphics[width=\textwidth]{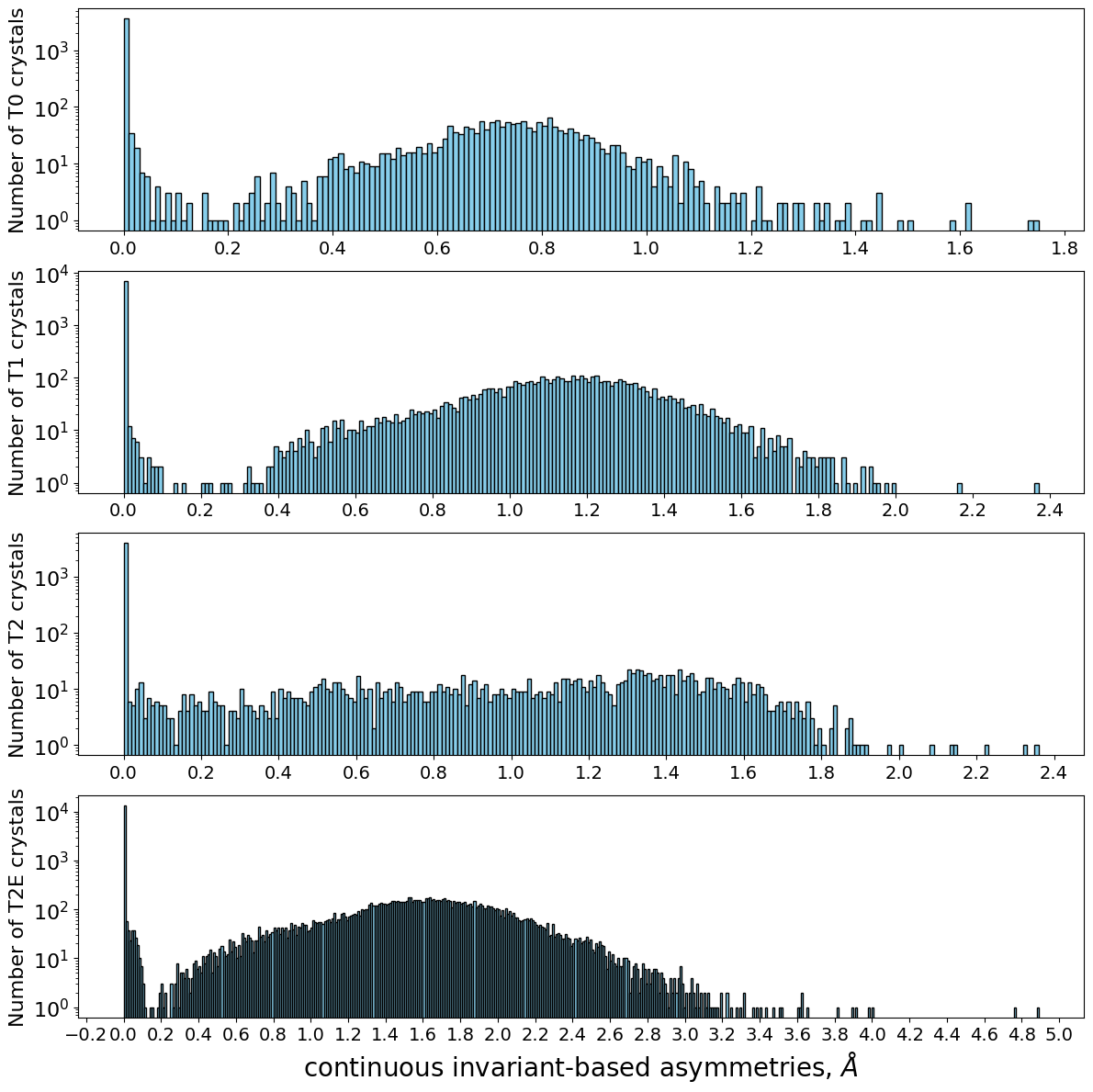}
\caption{The histograms of $\CIA_\infty$ for simulated crystals represented by 3 base
points at `ends' of molecules in Fig.~\ref{fig:T-molecules}.
\textbf{Row 1}: T0.
\textbf{Row 2}: T1.
\textbf{Row 3}: T2.
\textbf{Row 4}: T2E. 
}
\label{fig:T_CIA_inf}
\end{figure}

\begin{figure}
  \centering
  \includegraphics[width=\textwidth]{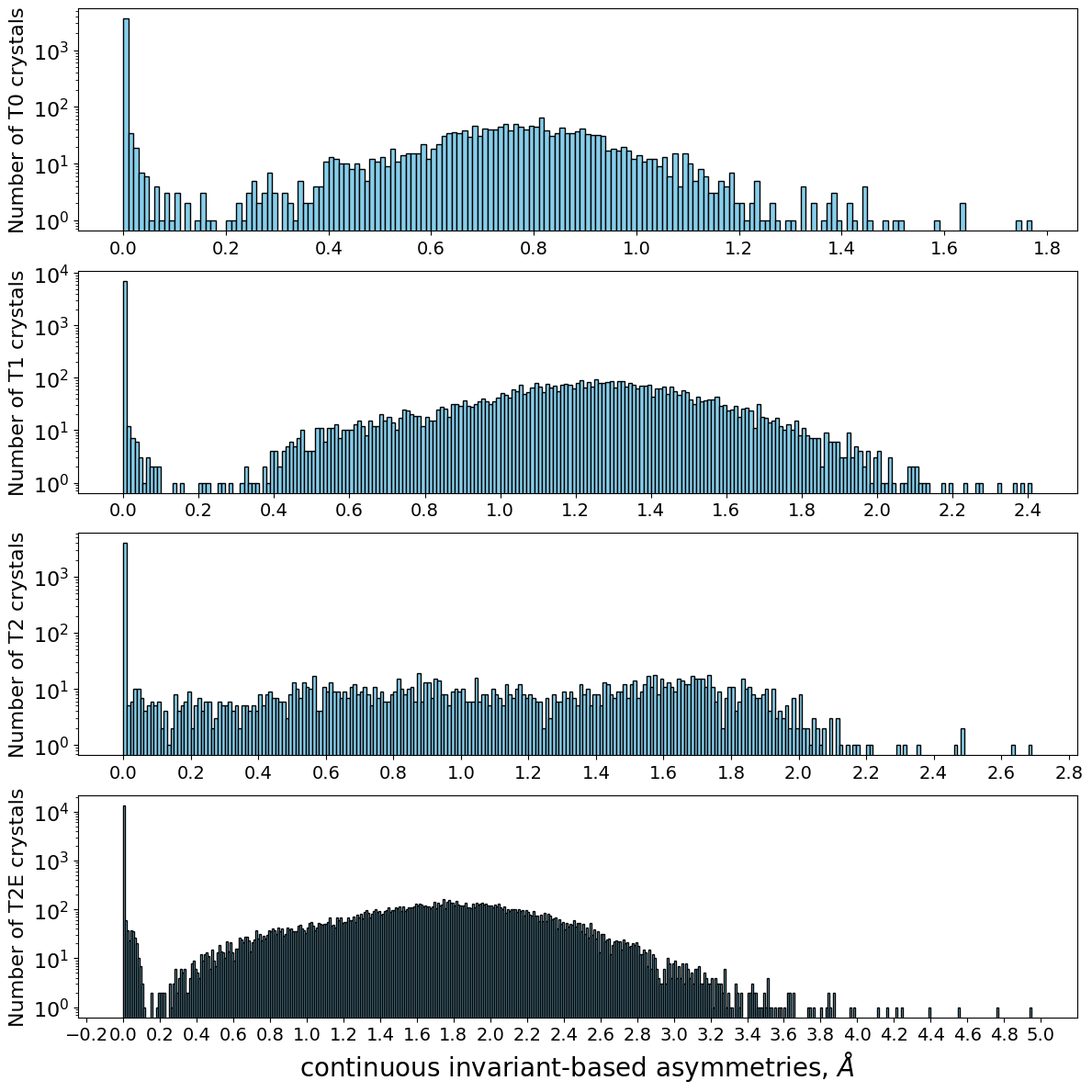}
\caption{The histograms of $\aCIA_\infty$ for simulated crystals represented by 3
base points at `ends' of molecules in Fig.~\ref{fig:T-molecules}.
\textbf{Row 1}: T0.
\textbf{Row 2}: T1.
\textbf{Row 3}: T2.
\textbf{Row 4}: T2E. 
}
\label{fig:T_aCIA_inf}
\end{figure}

\begin{table}
\label{tab:T-crystals_CIAs}
\centering
\caption{Statistics of $\aCIA,\CIA_\infty,\aCIA_\infty$ for the four CSP datasets
from \cite{pulido2017functional}.
The last rows contain Person correlations $r(x,y)$ between energy, density, and new
$\CIA$s.}
\medskip

\begin{tabular}{l|cccc}
CSP datasets & T0 crystals & T1 crystals & T2 crystals & T2E crystals \\ \hline
maximum $\CIA_\infty$, $\angstrom$ & 1.748 & 0.902 & 2.352 & 4.882 \\ 
$r(\text{energy}, \aCIA)$ & -0.393 & -0.196 & +0.035 & -0.020 \\ 
$r(\text{energy},\CIA_\infty)$ & $-0.398$ & $-0.196$ & $+0.016$ & $-0.019$ \\
$r(\text{energy}, \aCIA_\infty)$ & -0.399 & -0.186 & +0.032 & -0.014 \\
$r(\text{density}, \aCIA)$ & +0.315 & +0.144 & +0.036 & -0.021 \\
$r(\text{density},\CIA_\infty)$ & $+0.322$ & $+0.133$ & $+0.037$ & $-0.033$ \\
$r(\text{density}, \aCIA_\infty)$ & +0.323 & +0.131 & +0.032 & -0.022
\end{tabular}
\end{table}

\begin{figure}
  \centering
  \includegraphics[width=\textwidth]{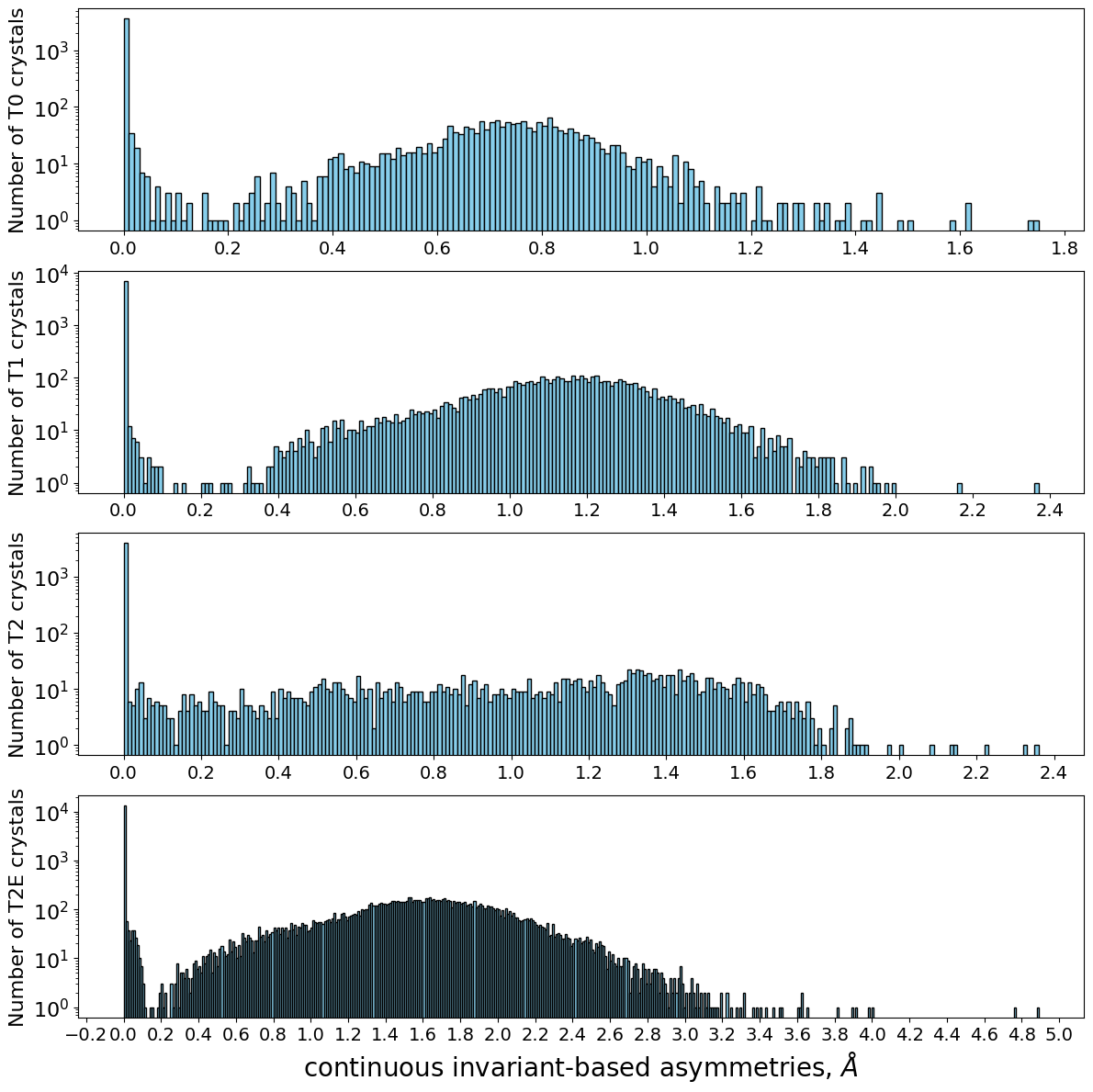}
\caption{The histograms of $\CIA_\infty$ for simulated crystals represented by 3 base
points at `ends' of molecules in Fig.~\ref{fig:T-molecules}.
\textbf{Row 1}: T0.
\textbf{Row 2}: T1.
\textbf{Row 3}: T2.
\textbf{Row 4}: T2E. 
}
  \label{fig:T_CIA_inf}
\end{figure}

\begin{figure}
    \centering
    \includegraphics[width=\textwidth]{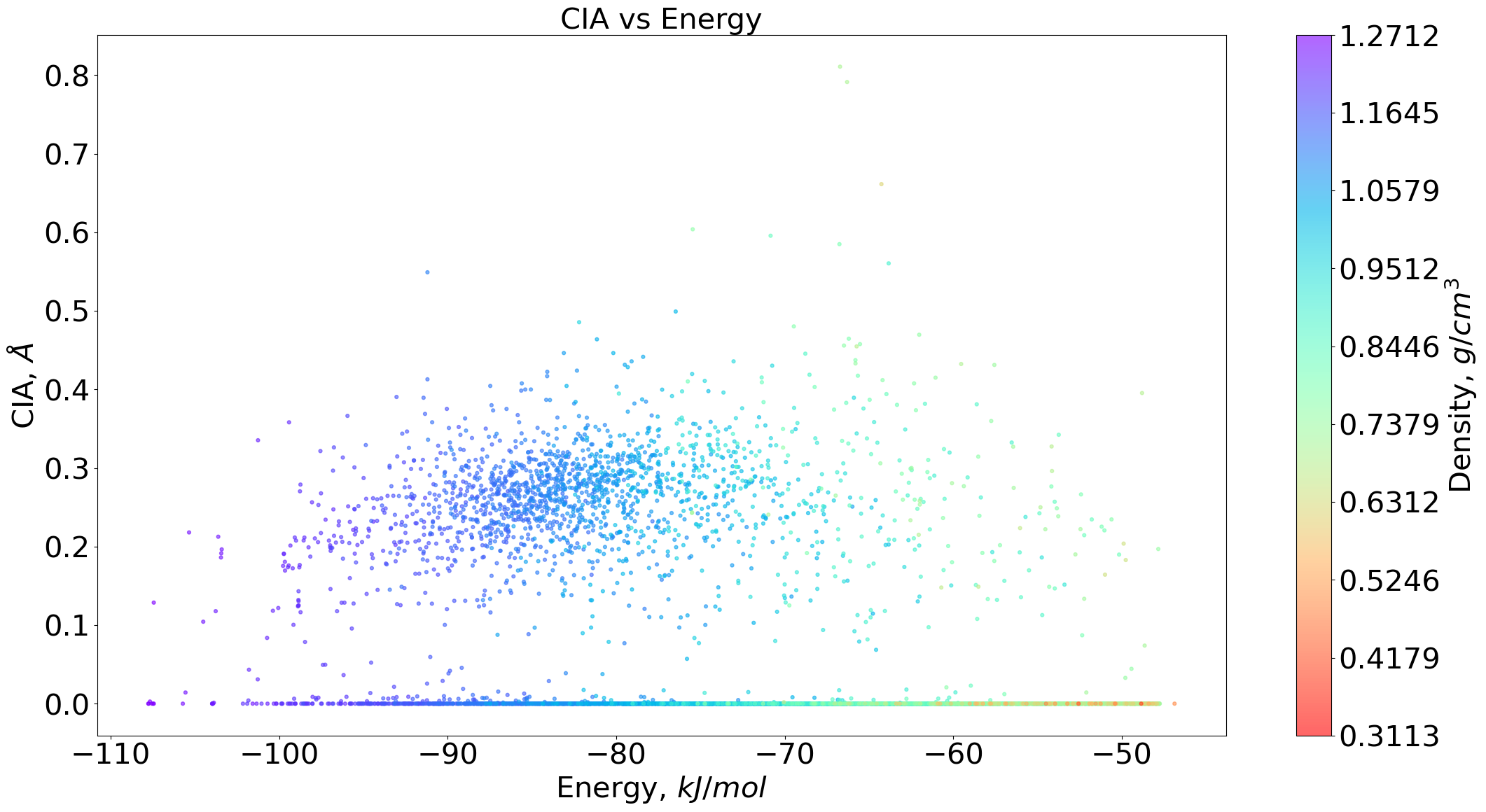}
\caption{CIA vs energy plot for simulated T0 crystals, coloured by their density}    \label{fig:T0_CIA_vs_energy}
\end{figure}

\begin{figure}
    \centering
    \includegraphics[width=\textwidth]{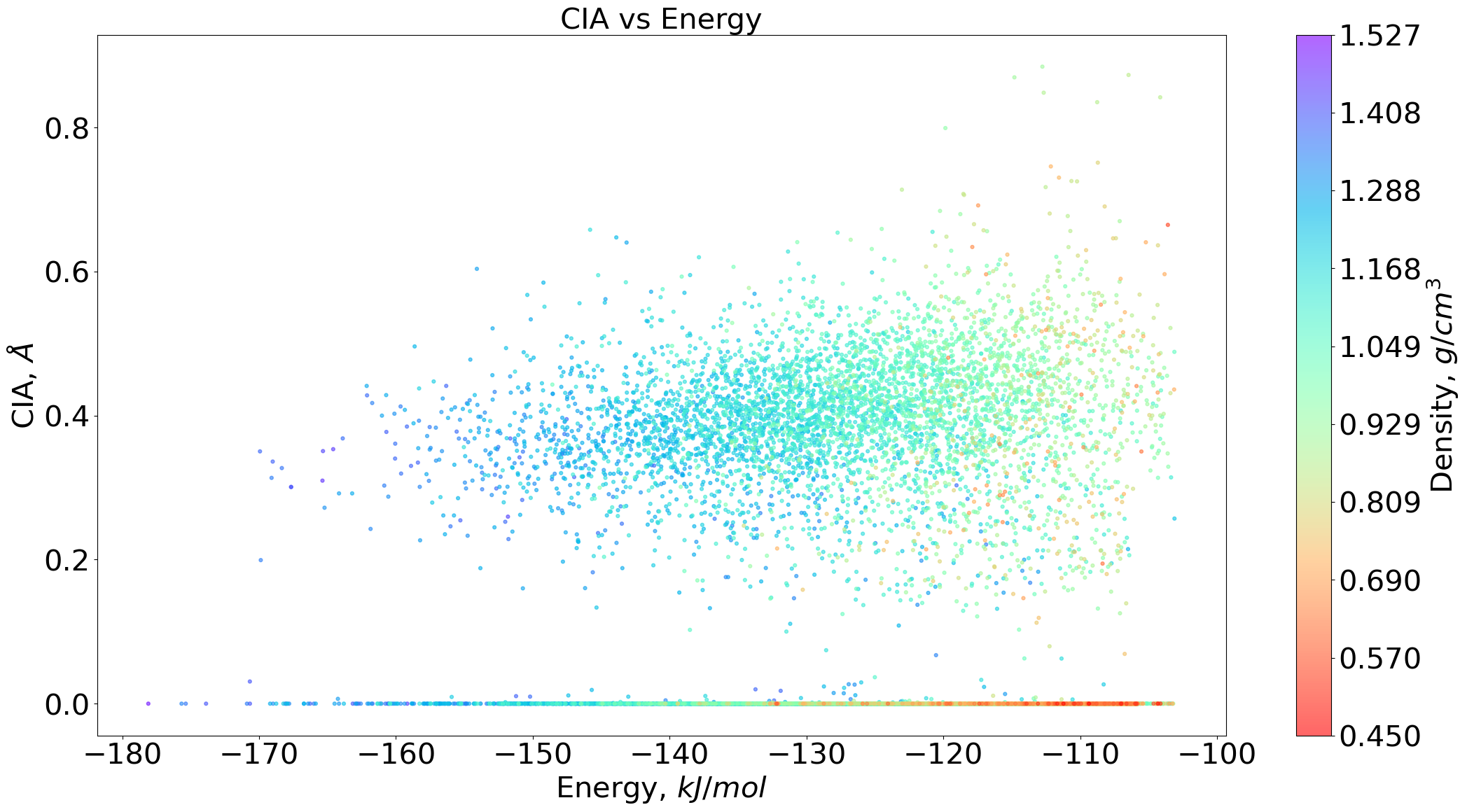}
\caption{CIA vs energy plot for simulated T1 crystals, coloured by their density}    \label{fig:T1_CIA_vs_energy}
\end{figure}

\begin{figure}
    \centering
    \includegraphics[width=\textwidth]{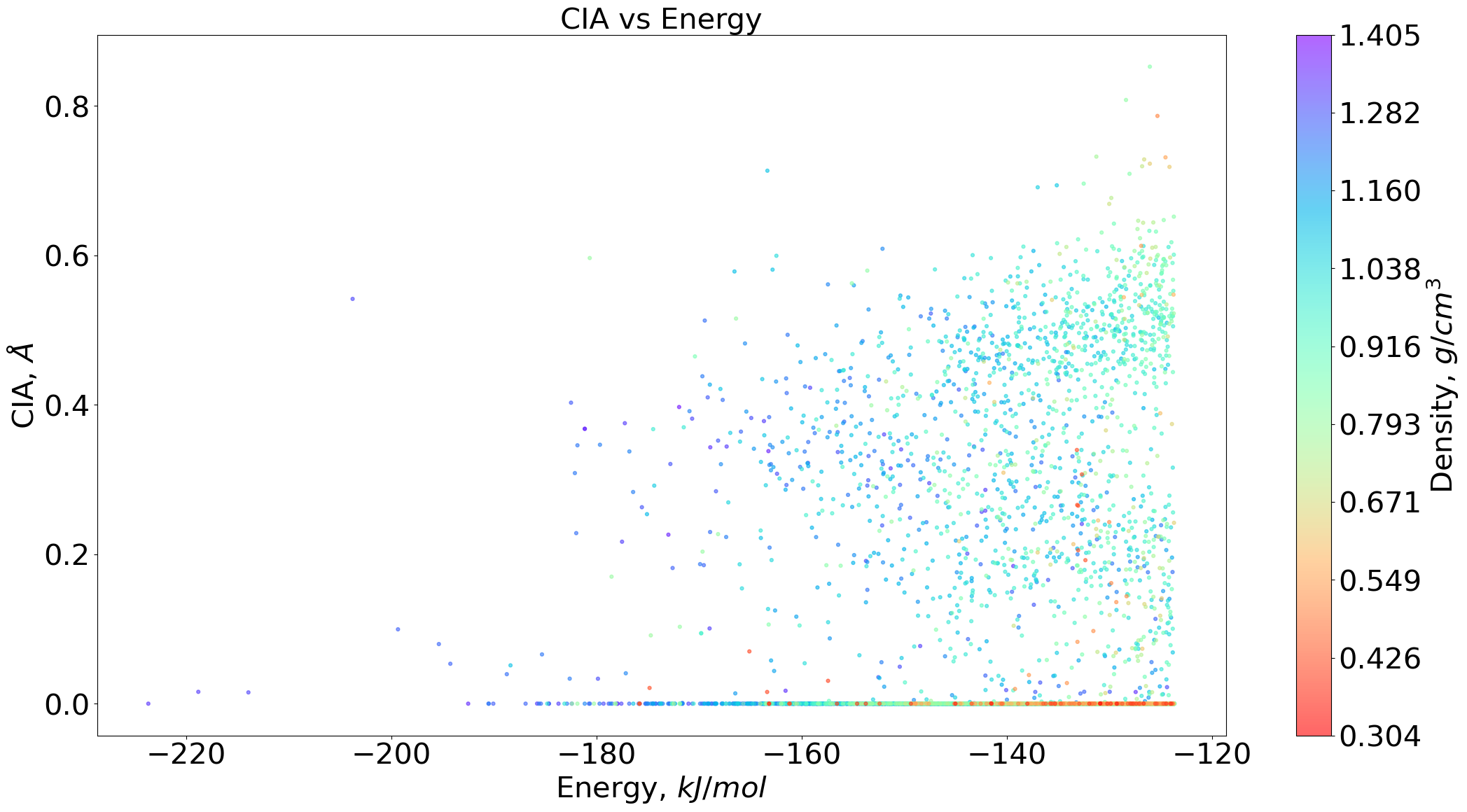}
\caption{CIA vs energy plot for simulated T2 crystals, coloured by their density}    \label{fig:T2_CIA_vs_energy}
\end{figure}

\begin{figure}
    \centering
    \includegraphics[width=\textwidth]{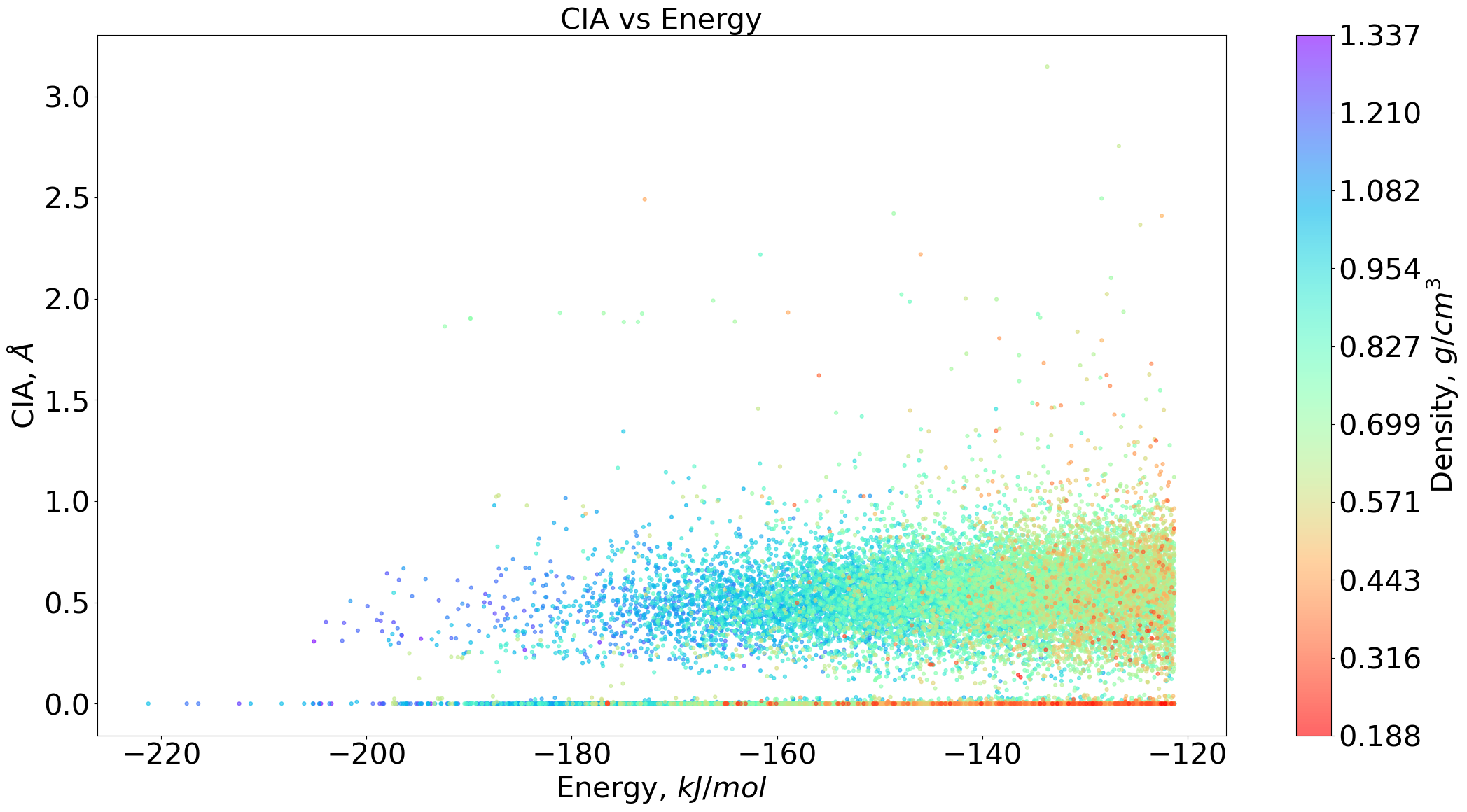}
\caption{CIA vs energy plot for simulated T2E crystals, coloured by their density}    \label{fig:T2E_CIA_vs_energy}
\end{figure}

Figures~\ref{fig:T0_CIA_vs_energy}, ~\ref{fig:T1_CIA_vs_energy},
~\ref{fig:T2_CIA_vs_energy}, and ~\ref{fig:T2E_CIA_vs_energy}
plot the CIA ($\text{\AA}$) and CIA$_\infty$ ($\text{\AA}$) vs lattice energy
(kJ/mol) for T0, T1, T2, and T2E simulated crystals, respectively.
Similar to Figures~\ref{fig:T0_CIA_vs_density}-\ref{fig:T2E_CIA_vs_density}, the
crystals are coloured with different density ($\text{g/cm}^3$) values.
The colour bar for density is shown at the right side of each plot. Similar to the
figures~\ref{fig:T0_CIA_vs_density}-~\ref{fig:T2E_CIA_vs_density}, the colour
gradient in these figures also indicate most of the crystals having lower energy,
i.e., are stable, are often associated with higher density values compared to the
crystals having higher energy which are mostly associated with lower density. While
symmetric crystals exist across the full range of energies and densities, crystals
with higher energies have larger asymmetries. High-density structures are observed
with both zero and non-zero asymmetry, indicating that density alone does not
determine symmetry. 
and higher asymmetry.
-178.09 $\text{kJ/mol}$ for T1, -123.7 to -223.7 $\text{kJ/mol}$ for T2, and -121.26
to -221.26 $\text{kJ/mol}$ for T2E predicted structures.

experimental structure SEMFAU. The CIA ($\text{\AA}$) and CIA$_\infty$ ($\text{\AA}$)
is zero for the structure, indicating symmetry. Since the structure is nano-porous
and has a low-density, the PPC of the structure is comparatively higher high average
radius of balls 'packed' inside the unit cell. $\text{ADA}_1$ ($\text{\AA}$) and
$\text{ADA}_2$ ($\text{\AA}$) representing how much $\text{AMD}_1$ and $\text{AMD}_2$
from the closest and second closest neighbour, respectively, deviate from the PPC,
scaled for number of neighbours.

\section{Detailed proofs of mathematical results}
\label{sec:proofs}

\begin{proof}[Proof of Lemma~\ref{lem:CIA_invariance}]
We prove the invariance of $\CIA(S)$.
The proof for other versions is almost identical. 
Since $\PDA(S;k)$ consists of inter-point distances, which are invariant under any
isometry, $\CIA(S)$ is also invariant.
If a unit cell $U$ of $S$ is transformed to another cell $U'$ by a matrix from
$\SL(\Z,n)$, there is a 1-1 correspondence between all points in the original motif
$S\cap U$ and the new motif $S\cap U'$ that respects all distances to neighbours.
Then $\PDA(S;k)$ and hence $\CIA(S)$ remain invariant.
If a unit cell $U$ is scaled up by an integer factor $c$, the original motif $M=S\cap
U$ transforms into the $c$-times larger motif $M_c$ containing $c$ isometric copies
of $M$.
The scaled-up asymmetric unit contains $c$ times more blocks $B_1,\cdots,B_{cG}$,
which can be considered as $c$ copies of the original blocks.
For each fixed $i=1,\dots,m$, the matrix of pairwise distances $\EMD$s between $cG$
blocks consists of $c\times c$ copies of the original matrix $G\times G$ of
distances.
The distances to the farthest units $d_{ij}=\max\limits_{j=1,\dots,cG} \EMD(B_i,B_j)$
are obtained by concatenating $c$ copies of the original vector $(d_{i1},\dots,
d_{iG})$.
Then the maximum and average values for each vector remain the same, so $\CIA(S)$ is
invariant.
\end{proof}

\begin{proof}[Proof of Lemma~\ref{lem:CIA_inequalities}]
The inequalities $\CIA(S)\leq\CIA_\infty(S)$ and $\aCIA(S)\leq\aCIA_\infty(S)$ hold,
because the RMS distance $d$ is bounded from above by the Chebyshev distance
$d_\infty$.
The inequality $\CIA(S)\leq\aCIA$ holds, because $\CIA(S)=\min\limits_{i=1,\dots,G}
d_i\leq \dfrac{1}{G}\sum\limits_{i=1}^G d_i=\aCIA(S)$.
To prove the inequality $\aCIA(S)\leq 2\CIA(S)$, let $B_i$ minimise
$d_i=\max\limits_{j=1,\dots,l} \EMD(B_i,B_j)=\CIA(S)$.
For $j,k=1,\dots,G$, the triangle 
inequality 
$$\EMD(B_k,B_j)\leq \EMD(B_k,B_i)+\EMD(B_i,B_j)\leq 2d_i=2\CIA(S)$$ implies that
$d_k=\max\limits_{j=1,\dots,G} \EMD(B_k,B_j)\leq 2\CIA(S)$ for each $k=1,\dots,G$.
Then
$\aCIA(S)=\dfrac{1}{G}\sum\limits_{k=1}^G d_k\leq 2\CIA(S)$.
\end{proof}

\begin{proof}[Proof of Theorem~\ref{thm:CIA}] 
By Lemma~4.1 in \cite{edelsbrunner2021density},
if periodic point sets $S,Q\subset\R^n$ are related by an $\ep$-perturbation, then
$S,Q$ have a common lattice.
Since $\CIA$ is invariant under changes of a unit cell by
Lemma~\ref{lem:CIA_invariance}, we can assume that $S,Q$ have the same number $m$ of
points in a common unit cell and equal Point Packing Coefficients
$\PPC(S)=\PPC(Q)$ in Definition~\ref{dfn:PDA}.
By Lemma SM3.4 in \cite{widdowson2026pointwise}, all corresponding elements of
$\PDD(S;k),\PDD(Q;k)$ differ by at most $2\ep$, which generalises the basic fact that
perturbing any two points up to $\ep$ changes the distance between them up to $2\ep$
by the triangle inequality.
The same upper bound of $2\ep$ holds for differences between all corresponding
elements of $\PDA(S;k),\PDA(Q;k)$ in Definition~\ref{dfn:PDA} due to
$\PPC(S)=\PPC(Q)$.
For both ground distances (Chebyshev and Root Mean Square) between rows of $k$
distances, the upper bound of $2\ep$ between corresponding distances $|b_i-c_i|\leq
2\ep$, $i=1\dots,k$, guarantees the same upper bound for
$d_\infty=\max\limits_{i=1,\dots,k}|b_i-c_i|\leq 2\ep$ and
$d=\sqrt{\dfrac{1}{k}\sum\limits_{i=1}^k (b_i-c_i)^2}\leq
\sqrt{\dfrac{1}{k}\sum\limits_{i=1}^k (2\ep)^2}=2\ep$.
\smallskip
 
Let $B_1,\dots,B_G$ be all geometric blocks in an asymmetric unit of $S$.
Denote by $C_1,\dots,C_{G}$ the corresponding blocks in an asymmetric unit of $Q$ so
that each $C_i$ is an $\ep$-perturbation of $B_i$ for $i=1,\dots,G$.
By the argument above, all $m_i$ corresponding points of $B_i$ and $C_i$ have
$2\ep$-close rows in $\PDA(S;k)$ and $\PDA(Q;k)$, respectively, for $i=1,\dots,G$.
Then $d(R_j(B_i),R_j(C_i))\leq 2\ep$ for $j=1,\dots,m_i$, where the ground distance
$d$ is Chebyshev or $\RMS$.
In the notations of Definition~\ref{dfn:EMD}, if we set $f_{jj}=\dfrac{1}{m_i}$ for
$j=1,\dots,m_i$, else $0$, then $\EMD(B_i,C_i)\leq 2\ep$.
The triangle inequality implies that $$\EMD(B_i,B_j)\leq
\EMD(B_i,C_i)+\EMD(C_i,C_j)+\EMD(C_j,B_j)|\leq \EMD(C_i,C_j)+4\ep.$$
Swapping the $B$-blocks and $C$-blocks, we similarly get $\EMD(C_i,C_j)\leq
\EMD(B_i,B_j)+4\ep$ and $|\EMD(B_i,B_j)-\EMD(C_i,C_j)|\leq 4\ep$, so the corresponding elements of the
matrix of $\EMD$s differ by at most $4\ep$.
Then the maximum distances $d_{i}$ in Definition~\ref{dfn:CIA} and hence the minima and averages over $j=1,\dots,G$ differ by at most $4\ep$.
\end{proof}

\end{document}